\documentclass{amsart}
\newtheorem{theorem}{Theorem}[subsection]
\newtheorem{lemma}[theorem]{Lemma}
\newtheorem{proposition}[theorem]{Proposition}
\newtheorem{corollary}[theorem]{Corollary} 
\theoremstyle{definition}  
\newtheorem{definition}[theorem]{Definition}
\newtheorem{example}[theorem]{Example}

\theoremstyle{remark}
\newtheorem{remark}[theorem]{Remark}

\newcommand{\id}{\mbox{id}} 
\newcommand{\Image}{\mbox{Image}}

\newcommand{\Map}{\mbox{Map}}
\newcommand{\End}{\mbox{End}} 
\newcommand{\Hom}{\mbox{Hom}} 
\newcommand{\Ad}{\mbox{Ad}\,}
\newcommand{\Rep}{\mbox{Rep}}
 
\renewcommand{\span}{\mbox{span}} 
    
\newcommand{\eps}{\varepsilon}
\newcommand{\bTheta}{{\bar \Theta}}
\newcommand{\mR}{{\mathcal R}} 
\newcommand{\bmR}{{\bar \mR}}

\newcommand{\bF}{{\bar F}}
\newcommand{\bJ}{{\bar J}}
\newcommand{\tJ}{{\mathcal J}}
\newcommand{\bK}{{\bar K}}
\newcommand{\tLambda}{\tilde{\Lambda}}
\newcommand{\g}{\mathfrak{g}}

\newcommand{\la}{{<}}  
\newcommand{\ra}{{>}}  

\newcommand{\1}{_{(1)}} 
\newcommand{\2}{_{(2)}} 
\newcommand{\3}{_{(3)}} 
 
\newcommand{\I}{^{(1)}} 
\newcommand{\II}{^{(2)}}

\begin{document}

\title{Dynamical quantum groups at roots of $1$}  

\author{Pavel Etingof}
\address{M.I.T., Department of Mathematics, 
Building 2, Room 2-165,
Cambridge, MA 02139-4307, USA}
\email{etingof@math.mit.edu}

\author{Dmitri Nikshych} 
\address{U.C.L.A., Department of Mathematics,  
405 Hilgard Avenue, Los Angeles, CA 90095-1555, USA}
\email{nikshych@math.ucla.edu}

\date{March 27, 2000}
\maketitle  

\begin{section}
{Introduction}
\label{intro}
 
{\bf 1.} The notion of a dynamical quantum group was first suggested by Felder
\cite{Fe} in 1994. Namely, Felder considered the quantum dynamical Yang-Baxter 
equation (also known as the Gervais-Neveu equation), 
which is a generalization of the usual quantum Yang-Baxter 
equation, and used the Faddeev-Reshetikhin-Takhtajan method to associate 
to a solution $R$ of this equation a certain algebra -- the dynamical
quantum group $F_R$. Felder also considered representations of $F_R$, 
and showed that although $F_R$ is not a Hopf algebra, its representations 
form a tensor category. 

In 1991, Babelon \cite{Ba} generalized Drinfeld's twisting method 
to the dynamical case, introducing the notion of a dynamical twist
(see also \cite{BBB}). 
Given a dynamical twist in a quasitriangular Hopf algebra $U$, 
one can define a solution of the dynamical Yang-Baxter 
equation acting on the tensor square of any representation of $U$. 

In 1997, it was shown independently in \cite{ABRR}, \cite{JKOS}, \cite{EV1}
that one can naturally construct a dynamical twist in the universal 
enveloping algebra of any simple Lie algebra or its $q$-deformation. 
Using the method of \cite{BBB}, one can obtain from these twists the 
solutions of the quantum dynamical Yang-Baxter equation from \cite{Fe}.  

At the same time, it was shown in \cite{EV1} that 
to any dynamical twist $J$, one can associate 
a Hopf algebroid $F(J)$. In the cases of \cite{Fe}, 
this Hopf algebroid coincides with $F_R$ as an algebra. 
In particular (as was shown already in 
\cite{EV2}), $F_R$ has a structure of a Hopf algebroid.
This explains the existence of  
a tensor product on the category of representations of $F_R$. 

In 1999, P.~Xu \cite{Xu1} associated to a dynamical twist 
$J$  on a Hopf algebra $U$,
another Hopf algebroid $U(J)$, 
obtained by twisting $U$ by means of $J$.
This Hopf algebroid $U(J)$ is 
closely related to the quasi-Hopf algebra associated to $(U,J)$ in 
\cite{BBB}, \cite{JKOS} (see \cite{Xu2}). P.Xu suggested that 
$U(J)$ should be dual, in an appropriate sense, to $F(J)$
(this is motivated by the duality of their classical limits). 
However, it is not very convenient 
to formulate such a statement precisely, because of difficulties with 
the notion of a dual Hopf algebroid. 

Moreover, it was shown in \cite{EV1} that for dynamical twists $J$ 
constructed in \cite{ABRR}, \cite{JKOS}, \cite{EV1}, the category 
of representations of $F(J)$ is essentially the same 
(as a tensor category) as that of $U(\g)$ or 
$U_q(\g)$. Thus, it is essentially the same as that of $U(J)$, 
which suggests that $U(J)$ and $F(J)$ should be not only 
dual to each other but also isomorphic. In other words, $U(J)$ and  $F(J)$ 
should be selfdual. Note that this would be a fundamentally new property, 
not satisfied by the usual Drinfeld-Jimbo quantum groups. 

{\bf 2.} This paper has two main goals: to make the above picture precise, 
and to generalize the theory of dynamical quantum groups to 
the case when the quantum parameter $q$ is a root of unity. 

Our first step is that we replace the notion of a Hopf algebroid 
with the recently introduced notion of a weak Hopf algebra 
\cite{BNSz}, \cite{BSz}. Roughly speaking, a weak Hopf algebra 
is an algebra and a coalgebra such that $\Delta$ is a homomorphism of algebras 
but is allowed to map $1$ to some idempotent not equal to $1\otimes 1$
(see Section~\ref{1} for a precise definition). 
Every weak Hopf algebra has a natural structure of a Hopf algebroid, 
but not vice versa. However, it turns out that dynamical quantum groups 
(at roots of unity) are a nice class of Hopf algebroids which do come from 
weak Hopf algebras. Moreover, regarding dynamical quantum groups 
as weak Hopf algebras rather than Hopf algebroids is convenient 
for two reasons: first, the definition of a weak Hopf algebra is much simpler, 
and second, it is naturally self-dual. 
 
Our main results can be summarized as follows. 

1. We generalize Drinfeld's twisting theory to weak Hopf algebras
(Section~\ref{2}). 
In particular, we show that twisting of a quasitriangular weak Hopf algebra
(defined in \cite{NVT}) gives another quasitriangular weak Hopf algebra. 

2. For every dynamical twist $J:\mathbb{T}\to U^{\otimes 2}$ 
of a Hopf algebra $U$, we define two weak Hopf algebras $H$ and $D$, 
the first by analogy with the construction of \cite{Xu1}, 
and the second by analogy with the construction of 
\cite{EV1} (they correspond to the Hopf algebroids
of \cite{Xu1}, \cite{EV1}; see Section~\ref{3}). 
We show that $H$ is isomorphic to $D^*$ with 
opposite multiplication. We consider the special case 
when $U$ is quasitriangular, and analyze the homomorphism 
$H^{*op}=D\to H$ defined by the quasitriangular structure on $H$. 
We give a criterion on when this homomorphism is an isomorphism. 

3. We take $U$ to be quantum group $U_q(\g)$ (for a finite dimensional 
simply laced Lie algebra $\g$), when $q$ is a root of unity
(more precisely, the finite dimensional version considered by Lusztig
\cite{L}; see Section~\ref{4}). 
We show that the known methods of producing dynamical twists 
for generic $q$ (\cite{ABRR}, \cite{ESS}) can also be used to produce 
dynamical twists when $q$ is a root of unity. In particular, 
for every generalized Belavin-Drinfeld triple $(\Gamma_1,\Gamma_2,T)$
for $\g$, we construct (following \cite{ESS})  a family of dynamical twists 
for $U_q(\g)$ which depends on $|\Gamma_1|$ parameters. With 
appropriate modifications this construction can be carried out in
the non-simply laced case as well.
 
4. We show that if $T$ is an automorphism of the Dynkin
diagram of $\g$ (in particular, if $T=1$) then the 
weak Hopf algebras $H$ and $D$ associated
to the corresponding twists are isomorphic 
(via the $R$-matrix of $H$). In particular, $D$ is self-dual 
(isomorphic to $D^{*,op}$), as was expected (for $T=1$) 
in the case of generic $q$. In particular, this implies that 
in this case the 
category of representations of the algebra $D=D_T$ corresponding to $T$ 
is equivalent (as a tensor category) to $\Rep(U)$, and thus to 
$\Rep(D_{T'})$ for any other automorphism $T'$. 
In particular, $\Rep(D_T)$ is equivalent to $\Rep(D_1)$.

Note that an analog of the latter result (when $q$ is generic,
$\g$ is the affine Lie algebra $\widehat{sl}_n$, and $T$ 
is the rotation of the Dynkin diagram of $\g$ by the angle 
$2\pi k/n$, $(k,n)=1$) is proved in \cite{ES2}.  

{\bf Acknowledgements.} We thank Ping Xu and Philippe Roche
for useful discussions.  The first author was supported by the 
NSF grant DMS-9700477. The second author thanks UCLA for providing 
him a research assistantship and MIT for the warm hospitality during his visit.

\end{section}

\begin{section}
{Weak Hopf algebras and Hopf algebroids}
\label{1}

Throughout this paper we work over an algebraically
closed field $k$ and use Sweedler's notation for
comultiplication, writing $\Delta(h) = h\1 \otimes h\2$.

\subsection{Weak Hopf algebras}

\begin{definition}[\cite{BNSz}, \cite{BSz}]
\label{weak bialgebra}
A {\em weak bialgebra} is  a $k$-vector space $H$ that has structures 
of an algebra $(H,\,m,\,1)$ and a coalgebra $(H,\,\Delta,\,\eps)$ 
such that the following axioms hold :
\begin{enumerate}
\item[1.] $\Delta$ is a (not necessarily unit-preserving) homomorphism :
\begin{eqnarray}
\Delta(hg) &=& \Delta(h)\Delta(g);
\label{eqn: delta homo}
\end{eqnarray}
\item[2.] The unit and counit satisfy the identities :
\begin{eqnarray}
\eps(hgf) &=& \eps(hg\1)\eps(g\2f) =  \eps(hg\2)\eps(g\1f),
\label{eqn: eps m} \\
(\Delta \otimes \id) \Delta(1) &=&
(\Delta(1)\otimes 1)(1\otimes \Delta(1)) =   
(1\otimes \Delta(1))(\Delta(1)\otimes 1),
\label{eqn: delta 1}
\end{eqnarray}
\end{enumerate}
for all $h,g,f\in H$.
\end{definition}

\begin{definition}[\cite{BNSz}, \cite{BSz}]
\label{basic definition} 
A {\em weak Hopf algebra}  $H$ is a weak bialgebra
equipped with a linear map $S: H \to H$, called an {\em antipode}, 
satisfying the following axioms :
\begin{eqnarray}
m(\id \otimes S)\Delta(h) &=&(\eps\otimes\id)(\Delta(1)(h\otimes 1)),
\label{eqn: S epst} \\
m(S\otimes \id)\Delta(h) &=& (\id \otimes \eps)((1\otimes h)\Delta(1)),
\label{eqn: S epss} \\
S(h\1)h\2 S(h\3)&=& S(h),
\label{S id S}
\end{eqnarray}
for all $h\in H$.
\end{definition}

Here axioms (\ref{eqn: eps m}) and (\ref{eqn: delta 1}) 
of Definition~\ref{weak bialgebra} are 
analogous to the bialgebra axioms of $\eps$ being an
algebra homomorphism and $\Delta$ a unit preserving map,
axioms (\ref{eqn: S epst}) and (\ref{eqn: S epss}) 
of Definition~\ref{basic definition} generalize the properties
of the antipode with respect to the counit. Also, it is possible to show
that given (\ref{eqn: delta homo}) - (\ref{eqn: S epss}),
axiom (\ref{S id S}) is equivalent to $S$ being both
anti-algebra and anti-coalgebra map.

A {\em morphism} of weak Hopf algebras is a map between them
which is both an algebra and a coalgebra morphism commuting with
the antipode. 

The image of a morphism is clearly
a weak Hopf algebra. The tensor product of two weak Hopf algebras
is defined in an obvious way.

Below we summarize the basic properties of weak Hopf algebras,
see \cite{BNSz} for the proofs.

The antipode $S$ of a weak Hopf algebra $H$ is unique; 
if $H$ is finite-dimensional then it is bijective \cite{BNSz}.

The right-hand sides of the formulas (\ref{eqn: S epst}) 
and (\ref{eqn: S epss}) are called the {\em target} 
and {\em source counital maps} and denoted 
$\eps_t$, $\eps_s$
respectively :
\begin{eqnarray}
\eps_t(h) = (\eps\otimes\id)(\Delta(1)(h\otimes 1)),\\
\eps_s(h) = (\id \otimes \eps)((1\otimes h)\Delta(1)).
\end{eqnarray}
The counital maps $\eps_t$ and $\eps_s$ are idempotents in $\End_k(H)$,
and satisfy  relations $S\circ \eps_t = \eps_s \circ S$ and
$S\circ \eps_s = \eps_t \circ S$.

The main difference between weak and usual Hopf algebras
is that the images of the counital maps are not necessarily equal 
to $k$. They turn out to be  subalgebras of $H$ 
called {\em target} and {\em source counital subalgebras} 
or {\em bases} as they generalize the notion of a base
of a groupoid  (cf.\ examples below) :
\begin{eqnarray}
H_t &=& \{h\in H \mid \eps_t(h) =h \}
  =  \{ (\phi\otimes \id)\Delta(1) \mid \phi\in H^* \}, \\
H_s &=& \{h\in H \mid \eps_s(h) =h \} 
  =  \{ (\id\otimes \phi)\Delta(1) \mid  \phi\in H^* \}.
\end{eqnarray}  
The counital subalgebras commute and  the restriction of the antipode
gives an anti-isomorphism between $H_t$ and $H_s$.

Any morphism between weak Hopf algebras preserves counital
subalgebras, i.e., if $\Phi : H\to H'$ is a morphism then
its restrictions on the counital subalgebras are isomorphisms :
$\Phi|_{H_t}: H_t\cong H'_t$ and $\Phi|_{H_s}: H_s\cong H'_s$.

The algebra $H_t$ (resp.\ $H_s$) is separable (and, therefore,
semisimple \cite{P}) with the separability idempotent
$e_t = (S\otimes \id)\Delta(1)$ (resp.\ $e_s=(\id\otimes S)\Delta(1)$),
i.e., we have $m(e_t) =m(e_s) =1$ and
\begin{eqnarray*}
(z_1\otimes 1)e_t(z_2\otimes 1) &=& (1\otimes z_2)e_t(1\otimes z_1),
\qquad z_1, z_2 \in H_t,\\
(y_1\otimes 1)e_s(y_2\otimes 1) &=& (1\otimes y_2)e_s(1\otimes y_1),
\qquad y_1, y_2 \in H_s.
\end{eqnarray*}
As a consequence of this fact and $\Delta$ being a homomorphism,
we have the following useful identities :
\begin{eqnarray}
(1\otimes z_1)\Delta(h)(1\otimes z_2)&=&
(S(z_1)\otimes 1)\Delta(h)(S(z_2)\otimes 1),
\qquad z_1, z_2 \in H_t,
\label{eqn: z delta} \\
(y_1\otimes 1)\Delta(h)(y_2\otimes 1) &=& 
(1\otimes S(y_1))\Delta(h)(1\otimes S(y_2)),
\qquad y_1, y_2 \in H_s.
\label{eqn: y delta}   
\end{eqnarray}

Note that $H$ is an ordinary Hopf algebra if and only if
$\Delta(1)=1\otimes 1$ if and only if $\eps$ is a homomorphism 
if and only if $H_t=H_s =k$.

The dual vector space $H^*$ has a natural structure of a weak
Hopf algebra with the structure operations dual to those of $H$ :
\begin{eqnarray}
& & \la \phi\psi,\,h\ra = \la \phi\otimes\psi,\,\Delta(h) \ra, \\
& & \la \Delta(\phi),\,h \otimes g \ra =  \la \phi,\,hg \ra, \\
& & \la S(\phi),\,h\ra = \la \phi,\,S(h) \ra, 
\end{eqnarray}
for all $\phi,\psi \in H^*,\, h,g\in H$. The unit of $H^*$ is $\eps$
and counit is $\phi \mapsto \la\phi,\, 1\ra$.

One can check that if $S$ is invertible, then
the opposite algebra $H^{op}$ is a weak Hopf algebra 
with the same coalgebra structure and the antipode $S^{-1}$. Similarly,
the coopposite coalgebra $H^{cop}$  is a weak Hopf algebra
with the same  algebra structure and the antipode $S^{-1}$.

It was shown in \cite{NVT} that modules over any weak Hopf
algebra $H$ form a monoidal category, called the {\em
representation category} and denoted $\Rep(H)$ with the product of
the modules $V$ and $W$ being equal to $\Delta(1)(V\otimes W)$
and the unit object given by $H_t$ which is an $H$-module via
$h\cdot z =\eps_t(hz),\,h\in H,z\in H_t$.

\begin{example}
\label{examples}
Let $G$ be a {\em groupoid} over a finite base (i.e., a category with 
finitely many objects, such that each morphism is invertible) 
then the groupoid algebra $kG$  is generated by morphisms $g\in G$ with 
the unit $1 =\sum_{X}\, \id_X$, where the sum is taken over all objects
$X$ of $G$, and  the product of  two morphisms is equal to  their composition
if the latter is defined and $0$ otherwise. 
It becomes 
a weak Hopf algebra via :
\begin{equation}
\Delta(g) = g\otimes g,\quad \eps(g) =1,\quad S(g)=g^{-1},\quad g\in G.
\end{equation}
The counital maps are given by $\eps_t(g) =gg^{-1}=\id_{target(g)}$ and
$\eps_s(g) =g^{-1}g = \id_{source(g)}$.

If $G$ is finite then 
the dual weak Hopf algebra $(kG)^*$ is generated by idempotents
$p_g,\, g\in G$ such that $p_g p_h= \delta_{g,h}p_g$ and
\begin{equation}
\Delta(p_g) =\sum_{uv=g}\,p_u\otimes p_v,\quad 
\eps(p_g)= \delta_{g,gg^{-1}}= \delta_{g,g^{-1}g},
\quad S(p_g) =p_{g^{-1}}.
\end{equation}
\end{example}
 
It is known that any group action on a set gives rise to a finite
groupoid. Similarly, in the non-commutative situation, one
can associate a weak Hopf algebra with every
action of a usual Hopf algebra on a separable algebra, see \cite{NVT}
for details.

\subsection{Hopf algebroids}
The following notions were introduced in \cite{Lu} (see also \cite{Xu1}).

\begin{definition}
\label{total algebra}
An algebra $H$ is called a {\em total algebra} with a {\em base algebra} 
$R$ if there exist
a {\em target map} $\alpha :R \to H$ which is an algebra homomorphism and
a {\em source map} $\beta :R \to H$ which is an algebra anti-homomorphism,
such that the images $\alpha(R)$ and $\beta(R)$ commute in $H$, i.e.,
\begin{equation}
\alpha(a)\beta(b) = \beta(b)\alpha(a), \quad \forall a,b\in R.
\end{equation}
\end{definition}

If $H$ is a total algebra then
above  maps define a natural $R-R$ bimodule structure on $H$ via
$$
a \cdot h \cdot b = \alpha(a) \beta(b) h, \quad h\in H,\,a,b\in R.
$$
Note that the bimodule tensor products $H\otimes_R H,\,
H\otimes_R H\otimes_R H,\dots$ are $R-R$ bimodules in an obvious way.

\begin{definition}
\label{comultiplication}
A {\em comultiplication} on a total algebra $H$ is an $R-R$ bimodule map
$\Delta : H\to H\otimes_R H$ 
satisfying $\Delta(1) = 1 \otimes_R 1$,
\begin{equation}
(\Delta\otimes_R \id)\Delta = (\id \otimes_R \Delta)\Delta : 
H \to H\otimes_R H\otimes_R H,
\end{equation}
and compatible with the maps $\alpha$, $\beta$ and comultiplication
in the sense that
\begin{eqnarray}
\Delta(h) (\beta(a)\otimes 1  - 1\otimes \alpha(a)) &=& 0, 
\label{eqn: compatibility 1} \\
\Delta(hg) &=& \Delta(h) \Delta(g),  \quad h,g\in H,\, a\in R.
\label{eqn: compatibility 2}
\end{eqnarray}
\end{definition}

Note that the right-hand side of equation (\ref{eqn: compatibility 2})
is well defined in $H\otimes_R H$ because of condition
(\ref{eqn: compatibility 1}).

\begin{definition}
\label{counit}
A {\em counit} for $H$ is an $R-R$ bimodule map 
$\epsilon: H\to R$ (where $R$ is an $R-R$ bimodule via multiplication) 
such that $\epsilon(1_H) =1_R$ and 
\begin{equation}
(\epsilon \otimes \id) \Delta = (\id \otimes \epsilon) \Delta = \id, 
\end{equation}
where we identify $R\otimes_R H \cong H \otimes_R R \cong H$.
\end{definition}

\begin{definition}
\label{bialgebroid}
A {\em bialgebroid} is a total algebra $H$ that possesses a comultiplication
and counit.
\end{definition}

\begin{definition}
\label{antipode}
An {\em antipode} for a bialgebroid $H$ with a base $R$ is a map
$\tau : H\to H$ which is an algebra anti-homomorphism  
such that $\tau\circ \beta =\alpha$ and satisfies the following properties :
\begin{enumerate}
\item[1.] 
$m(\tau\otimes \id)\Delta = \beta \circ \eps\circ \tau;$
\item[2.]
There exists a linear map $\gamma: H \otimes_R H \to H\otimes H$
which is a right inverse for the natural projection 
$H\otimes H \to H \otimes_R H$ such that
$m(\id\otimes \tau)\gamma\Delta = \alpha\circ \eps$.
\end{enumerate}
\end{definition}

Note that $m(\tau\otimes \id)$ is well defined on $H \otimes_R H$
but $m(\id\otimes \tau)$ is not, this is why it is necessary
to assume the existence of $\gamma$ in Definition~\ref{antipode}. 

\begin{definition}
\label{hopf algebroid}
A {\em Hopf algebroid} $H$ is a bialgebroid with an antipode.
\end{definition}

A (base preserving) morphism between Hopf algebroids
$H=(H,R,\dots)$ and $H=(H',R',\dots)$ is a pair
$(\psi,\Psi)$, where $\psi: R \to R'$ is an algebra isomorphism
and $\Psi: H\to H'$ is an algebra homomorphism which is also
an $R-R$ bimodule map such that 
$(\Psi\otimes\Psi)\Delta = \Delta'\circ \Psi$ and
\begin{eqnarray*}
\alpha'\circ\psi &=& \Psi \circ \alpha, \qquad
\beta'\circ\psi = \Psi \circ \beta,\\
\eps'\circ\Psi &=& \psi\circ \eps, \qquad
\tau'\circ\Psi = \Psi\circ \tau.
\end{eqnarray*} 

\subsection{The Hopf algebroid corresponding to a weak
Hopf algebra}

It turns out that weak Hopf algebras form a proper subclass 
of Hopf algebroids.

\begin{proposition}
\label{WHA is an algebroid}
Any weak Hopf algebra $H$ (not necessarily finite-dimensional)
has a natural structure of a Hopf algebroid (i.e., this assignment
defines a functor).
\end{proposition}
\begin{proof}
Define the base algebra to be the target subalgebra of $H$,
i.e., $R=H_t$, and let $\alpha =\id_{H_t},\, \beta = S^{-1}|_{H_t}$.
Then, clearly, $\alpha(R)= H_t,\, \beta(R)= H_s$, so that
images of $\alpha$ and $\beta$ commute. 

The comultiplication $\Delta$ regarded as a
map to $H\otimes_R H$ is coassociative and compatible with the multiplication.
It is a bimodule map since 
\begin{equation}
\Delta(\alpha(a)h) = (\alpha(a)\otimes 1)\Delta(h), \qquad
\Delta(\beta(a)h) = (1\otimes \beta(a))\Delta(h),
\end{equation}
for all $h\in H$ and $a\in R$.
We also have the compatibility condition:
$$
\Delta(h)(\beta(a)\otimes 1  - 1\otimes \alpha(a)) 
= h\1S^{-1}(a) \otimes h\2 - h\1 \otimes h\2a = 0.
$$
Next, let $\epsilon = \eps_t$, then we have $\epsilon(1)=1$ and
$$ 
\epsilon(a\cdot h \cdot b) = \eps_t(a S^{-1}(b)h) = a \epsilon(h) b,
$$
for all $h\in H,\, a,b\in R$, i.e., $\epsilon$ is an $R-R$ module map, also
\begin{eqnarray*}
\epsilon(h\1)\cdot h\2 &=& \eps_t(h\1)h\2 = \eps(1\1h\1)1\2h\2 = h,\\
h\1 \cdot \epsilon(h\2)&=& S^{-1}(\epsilon(h\2))h\1  = 
1\1 h\1 \eps(S(1\2) h\2)= h,
\end{eqnarray*}
where we used the antipode axioms of a weak Hopf algebra. Thus,
$\epsilon $ satisfies the counit axiom.

The antipode $\tau =S$ is an algebra anti-homomorphism satisfying
$\tau\circ \beta = \id_{H_t} =\alpha$. The section
$\gamma : H\otimes_R H \to H\otimes H$ is given by
$$
\gamma(h\otimes_R g) = (h\cdot S(1\1)) \otimes (1\2\cdot g) 
= \Delta(1)(h\otimes g).
$$
Finally, we verify the antipode properties : 
\begin{eqnarray*}
m(\tau \otimes_R \id)\Delta &=&\eps_s = S^{-1}\circ \eps_t\circ S 
= \beta\circ\epsilon\circ\tau, \\
m(\id \otimes_R  \tau)\gamma\Delta &=& \eps_t =\alpha\circ \epsilon.
\end{eqnarray*}
Thus, $(H, H_t, \id_{H_t}, S^{-1}_{H_t},\Delta, \eps_t, S)$
is a Hopf algebroid.

If $\Psi : H\to H'$ is a morphism between weak Hopf algebras, 
then it is clear from our definitions that the pair 
$(\Psi|_{H_t},\, \Psi)$ gives a morphism between the corresponding 
Hopf algebroids.
\end{proof}

\begin{remark}
\label{converse is false}
The converse of Proposition~\ref{WHA is an algebroid} is false
even when $H$ is finite dimensional. Indeed, the base algebra
$R=H_t$ of $H$ is necessarily separable, on the other hand,
for any algebra $A$ the space $A\otimes A^{op}$ has a structure 
of a Hopf algebroid with the base $A$ (\cite{Lu}, Examples 3.1
and 4.4). In the case when $A$ is not separable, this Hopf
algebroid is not a weak Hopf algebra.
\end{remark}

\end{section}


\begin{section}
{Weak Hopf algebras coming from twisting}
\label{2}

\subsection{Twisting}
We describe the procedure of constructing new weak Hopf algebras
by twisting a comultiplication. Twisting of Hopf algebroids without
an antipode was developed in \cite{Xu1} and  a special case of
twisting of weak Hopf $*$-algebras was considered in \cite{NV}.

\begin{definition}
\label{twist}
A {\em twist} for a weak Hopf algebra $H$ is a pair ($\Theta,\,\bTheta$), with
\begin{equation}
\label{Theta and bTheta} 
\Theta\in \Delta(1)(H\otimes H), \quad
\bTheta\in (H\otimes H)\Delta(1), \quad
\mbox { and } \quad
\Theta\bTheta= \Delta(1)
\end{equation}
satisfying the following axioms :
\begin{eqnarray}
(\eps\otimes \id)\Theta = (\id\otimes \eps)\Theta 
&=& (\eps\otimes \id)\bTheta = (\id\otimes \eps)\bTheta =1,
\label{eqn: twist eps} \\
(\Delta\otimes\id)(\Theta) (\Theta\otimes 1) &=&
(\id\otimes \Delta)(\Theta) (1\otimes\Theta), 
\label{eqn: twist ++}\\
(\bTheta\otimes 1) (\Delta\otimes\id)(\bTheta) &=&
(1\otimes \bTheta) (\id\otimes \Delta)(\bTheta), 
\label{eqn: twist --}\\
(\Delta\otimes\id)(\bTheta) (\id\otimes \Delta)(\Theta)  &=&
(\Theta\otimes 1) (1\otimes \bTheta), 
\label{eqn: twist +-}\\
(\id\otimes \Delta)(\bTheta) (\Delta\otimes\id)(\Theta) &=&
(1\otimes\Theta) (\bTheta\otimes 1).
\label{eqn: twist -+}
\end{eqnarray}
\end{definition}
For ordinary Hopf algebras this notion coincides with the usual
notion of twist and each of the four conditions~(\ref{eqn: twist ++})
-- (\ref{eqn: twist -+}) implies the other 
three. But since $\Theta$ and $\bTheta$ are, in general,
not invertible we need to impose all of them.

The next Proposition extends Drinfeld's  twisting construction 
to the weak case.

\begin{proposition}
\label{twisting}
Let ($\Theta,\bTheta)$  be a twist for a weak Hopf algebra $H$. 
Then there is a weak Hopf algebra $H_\Theta$ having
the same algebra structure and counit as $H$ with
a comultiplication and antipode given by
\begin{equation}
\Delta_\Theta(h) = \bTheta\Delta(h)\Theta, \qquad
S_\Theta(h) = v^{-1}S(h)v,
\label{eqn: theta}
\end{equation}
for all $h\in H_\Theta$, where $v= m(S\otimes\id)\Theta$
is invertible in $H_\Theta$.
\end{proposition}
\begin{proof}
Clearly, $\Delta_\Theta$ is an algebra homomorphism.
Its coassociativity follows from axioms~(\ref{eqn: twist ++})
and  (\ref{eqn: twist --}).

Observe that the relations between $\Theta$ and $\eps$
yield the following identities (recall that $S$ is invertible
when restricted on the counital subalgebras, since the latter
are finite-dimensional) :
\begin{eqnarray}
\eps_t(\Theta\I )\Theta\II &=& 1, \qquad
S^{-1}\eps_t(\Theta\II)\Theta\I  =1, \\
\bTheta^{(1)} \eps_s(\bTheta^{(2)}) &=& 1, \qquad
\bTheta^{(2)} S^{-1}\eps_s(\bTheta^{(1)}) = 1. 
\end{eqnarray}
Here and in what follows we write $\Theta = \Theta\I \otimes \Theta\II$ and
$\bTheta = \bTheta^{(1)} \otimes \bTheta^{(2)}$.

Using these properties and the equation~(\ref{eqn: z delta})
we check that $\eps$ is still a counit for ($H, \Delta_\Theta$):
\begin{eqnarray*}
(\eps \otimes \id)\Delta_\Theta(h)
&=& \eps(\bTheta^{(1)}h\1 \Theta\I ) \bTheta^{(2)}h\2\Theta\II \\
&=& \bTheta^{(2)} S^{-1}\eps_s(\bTheta^{(1)}) h \eps_t(\Theta\I )\Theta\II=h,\\
(\id \otimes \eps)\Delta_\Theta(h)
&=& \bTheta^{(1)}h\1 \Theta\I \eps( \bTheta^{(2)}h\2\Theta\II ) \\
&=& \bTheta^{(1)} \eps_s(\bTheta^{(2)}) h S^{-1}\eps_t(\Theta\II)\Theta\I=h. 
\end{eqnarray*}
We proceed to verify the rest of the axioms
of a weak Hopf algebra, writing
$\Delta_\Theta(h) = {h_\Theta}\1 \otimes {h_\Theta}\2$ :
\begin{eqnarray*}
\eps(g{h_\Theta}\1 )\eps({h_\Theta}\2f)
&=& \eps(g\bTheta^{(1)}h\1 \Theta\I)\eps(\bTheta^{(2)}h\2\Theta\II f)\\
&=& \eps(g\bTheta^{(1)} \eps_s(\bTheta^{(2)}) h\1)
    \eps(h\2 \eps_t(\Theta\I )\Theta\II f) \\
&=& \eps(gh\1)\eps(h\2 f) = \eps(ghf), \\
\eps(g{h_\Theta}\2 )\eps({h_\Theta}\1 f)
&=& \eps(g\bTheta^{(2)}h\2\Theta\II)\eps(\bTheta^{(1)}h\1 \Theta\I f)\\
&=& \eps(g \bTheta^{(2)} S^{-1}\eps_s(\bTheta^{(1)}) h\2)
    \eps(h\1 S^{-1}\eps_t(\Theta\II)\Theta\I f) \\
&=& \eps(g h\2)\eps(h\1 f)  = \eps(ghf),
\end{eqnarray*}
for all $g,h,f \in H$. 

The axioms involving $\Delta_\Theta(1)$ are checked using
identities (\ref{eqn: twist +-}) and (\ref{eqn: twist -+})
of Definition~\ref{twist} :
\begin{eqnarray*}
(1\otimes \Delta_\Theta(1)) (\Delta_\Theta(1))\otimes 1)
&=& (1\otimes \bTheta\Theta) (\bTheta\Theta\otimes 1) \\
&=& (1\otimes \bTheta) (\id\otimes \Delta)\bTheta
    (\Delta\otimes\id)\Theta (\Theta\otimes 1) \\
&=& (\Delta_\Theta \otimes \id) \Delta_\Theta(1) \\
&=& (\bTheta\otimes 1)(\Delta\otimes\id)\bTheta
    (\id\otimes \Delta)\Theta (1\otimes\Theta) \\
&=& (\bTheta\Theta)\otimes 1)(1\otimes \bTheta\Theta)\\
&=& (\Delta_\Theta(1))\otimes 1)(1\otimes \Delta_\Theta(1)).
\end{eqnarray*}  
We define a new antipode by
\begin{equation}
S_\Theta(h)= \bTheta^{(1)} S(\bTheta^{(2)} ) S(h)S(\Theta\I) \Theta\II
\end{equation}
and compute (writing $\Theta',\, \bTheta'$ for
additional copies of $\Theta$ and $\bTheta$) :
\begin{eqnarray*}
m(\id \otimes S_\Theta)\Delta_\Theta(h)
&=& \bTheta^{(1)} h\1 \Theta\I \bTheta^{'(1)}S(\bTheta^{'(2)})S(\Theta\II) \\
& &    S(h\2) S(\bTheta^{(2)})S(\Theta^{'(1)})\Theta^{'(2)} \\
&=& \bTheta^{(1)}h\1 S(h\2) S(\bTheta^{(2)})S(\Theta^{'(1)})\Theta^{'(2)}\\
&=& \bTheta^{(1)} \eps_t(h) S(\Theta\I \bTheta^{2})\Theta\II \\
&=& \bTheta^{(1)}  {\Theta\I}\1 \eps_t(h) S({\Theta\I}\2)
     \eps_s(\bTheta^{2})\Theta\II \\
&=& \bTheta^{(1)} \eps_s(\bTheta^{2}) \eps_t(\Theta\I h) \Theta\II \\
&=& \eps_t(\Theta\I h) \Theta\II \\
&=& \eps(1\1 \Theta\I h) 1\2 \Theta\II \\
&=& \eps(\bTheta^{(1)} 1\1 \Theta\I h) \bTheta^{(2)} 1\2 \Theta\II \\
&=& (\eps\otimes\id)(\Delta_\Theta(1)(h\otimes 1)), 
\end{eqnarray*}
where we used axioms of $(\Theta, \bTheta)$ and properties of the
counital maps. Note that for $h=1$ we get the identity
$$
\bTheta^{(1)}S(\bTheta^{(2)}) \cdot S(\Theta^{(1)})\Theta^{(2)} =1, 
$$
i.e., $\bTheta^{(1)}S(\bTheta^{(2)})$ is the inverse of
$v= S(\Theta^{(1)})\Theta^{(2)}$. It can be proven similarly that 
$$
m(S_\Theta\otimes \id)\Delta_\Theta(h)=
(\id \otimes \eps)((1\otimes h)\Delta_\Theta(1)), \quad h\in H.
$$
Finally, let us write
$$
(\Delta_\Theta \otimes \id)\Delta_\Theta(h) =
\bTheta^{(1)} h\1 \Theta^{(1)} \otimes \bTheta^{(2)} h\2 \Theta^{(2)}
\otimes \bTheta^{(3)} h\3 \Theta^{(3)}
$$
and compute
\begin{eqnarray*}
\lefteqn{ m(m\otimes \id)(S_\Theta \otimes \id \otimes S_\Theta)
\Delta_\Theta(h) =} \\
&=& v^{-1} S(\Theta^{(1)})S(h\1)S(\bTheta^{(1)})S(\Theta^{'(1)})\Theta^{'(2)}
    \bTheta^{(2)} \\
& & h\2 \Theta^{(2)}\bTheta^{'(1)}S(\bTheta^{'(2)})
    S(\Theta^{(3)})S(h\3) S(\bTheta^{(3)}) v \\
&=& v^{-1}S(\Theta^{(1)})S(h\1) \eps_s(\bTheta^{(1)})h\2 \eps_t(\Theta^{(2)})
    S(h\2) S(\bTheta^{(2)}) v \\
&=& v^{-1}S(\Theta^{(1)}) \eps_s(\bTheta^{(1)} h\1) \eps_t(\Theta^{(2)})
    S(h\2) S(\bTheta^{(2)}) v \\
&=& v^{-1}  \eps_s(\bTheta^{(1)} h\1)S(h\2) S(\bTheta^{(2)}) v \\
&=& v^{-1} S(h\1) \eps_s(\bTheta^{(1)} )\eps_t(h\2) S(\bTheta^{(2)}) v
    = v^{-1} S(h)v = S_\Theta(h),
\end{eqnarray*}
whence $H_\Theta$ is a weak Hopf algebra.
\end{proof}

\begin{remark}
\label{twisted bases}
Note that it follows from the proof of Proposition~\ref{twisting}
that the counital maps of the twisted weak Hopf 
algebra $H_\Theta$ are given by
\begin{eqnarray}
(\eps_t)_\Theta(h) &=& \eps(\Theta\I h)\Theta\II, \\
(\eps_s)_\Theta(h) &=& \bTheta^{(1)}  \eps(h \bTheta^{(2)}),
\end{eqnarray}
so the counital subalgebras are also getting 
deformed in general.
\end{remark}

\begin{remark}
\label{gauge equivalence}
If $(\Theta,\bTheta)$ is a twist for $H$ and $x\in H$ is an invertible
element such that $\eps_t(x)=\eps_s(x)=1$ then $(\Theta^x,\bTheta^x)$, where
$$
\Theta^x =\Delta(x)^{-1}\Theta (x\otimes x) \quad \mbox{ and } \quad
\bTheta^x = (x^{-1}\otimes x^{-1})\bTheta \Delta(x),
$$
is also a twist for $H$. The twists $(\Theta,\bTheta)$ and
$(\Theta^x,\bTheta^x)$ are called {\em gauge equivalent} and
$x$ is called a {\em gauge transformation}.  Given such an $x$,
the map $h\mapsto x^{-1}hx$ is an isomorphism between weak Hopf
algebras $H_\Theta$ and $H_{\Theta^x}$.
\end{remark}

\begin{remark}
If $\Theta$ is a twist for a weak Hopf algebra $H$,
then it also defines a twist of the corresponding
Hopf algebroid constructed as in Proposition~\ref{twisting},
cf.\ \cite{Xu1}.
\end{remark}

\subsection{Quasitriangular  weak Hopf algebras}
The notion of a {\em quasitriangular  weak Hopf algebra}
was introduced and studied in \cite{NVT}. It is defined
as a triple ($H,\, \mR,\,\bmR$) where $H$ is
a weak Hopf algebra with
\begin{equation} 
\mR\in \Delta^{op}(1)(H\otimes H)\Delta(1), \qquad
\bmR\in \Delta(1)(H\otimes H)\Delta^{op}(1), 
\end{equation}
\begin{equation}
\mR\bmR = \Delta^{op}(1), \qquad \bmR\mR = \Delta(1), \qquad
\mbox{ and } \qquad \Delta^{op}(h)\mR = \mR\Delta(h),
\end{equation}
for all $h\in H$,
where $\Delta^{op}$ denotes the comultiplication opposite to $\Delta$, 
and such that $\mR$ obeys the following conditions :
\begin{equation}
(\id \otimes \Delta)\mR = \mR_{13}\mR_{12}, \qquad
(\Delta \otimes \id)\mR = \mR_{13}\mR_{23},
\end{equation}
where $\mR_{12} = \mR\otimes 1$ etc.\ as usual.

The existence of a quasitriangular structure $\mR$ on $H$
is equivalent to  $\Rep(H)$  being a braided  category,
and any quasitriangular structure $\mR$ is a solution of the 
quantum Yang-Baxter equation :
\begin{equation}
\mR_{12}\mR_{13}\mR_{23} = \mR_{23}\mR_{13}\mR_{12}.
\end{equation}
An example of a quasitriangular weak Hopf
algebra is given by the Drinfeld double of finite dimensional 
weak Hopf algebra \cite{NVT}. 

As in the case of ordinary Hopf algebras, a quasitriangular
structure $\mR$ on $H$ defines two homomorphisms of weak Hopf algebras :
\begin{eqnarray}
&\rho_1 :& H^* \ni \phi \mapsto (\id \otimes \phi)(\mR) \in H^{op}, 
\label{eqn: H* to H} \\
&\rho_2 :& H^* \ni \phi \mapsto  (\phi\otimes \id)(\mR) \in H^{cop},
\end{eqnarray}
in particular, when $\mR$ has a maximal ($=\dim H$) rank,
then it defines isomorphisms $H^* \cong H^{op}\cong H^{cop}$.

A twisting of a quasitriangular weak Hopf algebra is again
quasitriangular. Namely, if $(\Theta,\bTheta)$ is a twist 
and $(\mR,\bmR)$ is a quasitriangular structure  for $H$ then 
the quasitriangular structure for $H_\Theta$ is given by 
$(\bTheta_{21}\mR \Theta, \bTheta \bmR \Theta_{21})$.
The proof of this fact is exactly the same as for ordinary Hopf algebras. 
\end{section}


\begin{section}
{Weak Hopf algebras arising from dynamical twists}
\label{3}

We describe two methods of constructing weak Hopf algebras,
which are finite-dimensional modifications of constructions
proposed in \cite{Xu1} and \cite{EV1}. It turns out that resulting
weak Hopf algebras are dual to each other.

\subsection{Dynamical twists on Hopf algebras}
\label{section 3.1}
Dynamical twists  first appeared in the work of Babelon \cite{Ba},
see also \cite{BBB}.

Let $U$ be a Hopf algebra and $A=\Map(\mathbb{T},\, k)$
be a commutative and cocommutative Hopf algebra
of functions on a finite Abelian group $\mathbb{T}$
which is a Hopf subalgebra of $U$. Let $P_\mu,\, \mu\in \mathbb{T}$ 
be the minimal idempotents in $A$. We fix this notation for
the rest of the paper.

\begin{remark} 
In \cite{EV1} the role of the group  $\mathbb{T}$ is played by the
dual space of a Cartan subalgebra of a simple Lie algebra.
\end{remark}

\begin{definition}
\label{zero weight}
We say that an elment $x$ in $U^{\otimes n},\, n\geq 1$ 
has {\em zero weight} if $x$ commutes with 
$\Delta^n(a)$ for all $a\in A$, where $\Delta^n: A \to A^{\otimes n}$ 
is the iterated comultiplication.
\end{definition}

\begin{definition}
\label{dynamical twist}
An invertible, zero-weight  $U^{\otimes 2}$-valued function 
$J(\lambda)$ on $\mathbb{T}$ is called a {\em dynamical twist } 
for $U$ if it satisfies
the following functional equations :
\begin{eqnarray}
\label{dynamical equation}
(\Delta\otimes \id)J(\lambda) (J(\lambda+h^{(3)}) \otimes 1)
&=& (\id \otimes\Delta) J(\lambda) (1\otimes J(\lambda) ), \\
(\eps\otimes \id)J(\lambda) = (\id \otimes \eps)J(\lambda)  &=& 1.
\end{eqnarray}
Here an in what follows the notation $\lambda+ h^{(i)}$ means that
the argument $\lambda$ is shifted by the weight of the $i$-th component,
e.g., $J(\lambda+h^{(3)}) =\sum_\mu\, J(\lambda+\mu) \otimes P_\mu
\in U^{\otimes 2} \otimes A$.
\end{definition}

\begin{remark}
\label{dynamical gauge}
If $J(\lambda)$ is a dynamical twist for $U$ and $x(\lambda)$
is an invertible,  zero-weight  $U$-valued function on $\mathbb{T}$
such that $\eps(x(\lambda))\equiv 1$, then 
$$ 
J^x(\lambda) = \Delta(x(\lambda)^{-1})\,J(\lambda)\, 
                (x(\lambda+h^{(2)})\otimes x(\lambda))
$$
is also a dynamical twist for $U$, {\em gauge equivalent} to $J(\lambda)$.  
\end{remark}

Note that for every fixed $\lambda\in \mathbb{T}$ the element
$J(\lambda)\in U\otimes U$ does not have to be a twist for $U$
in the sense of Drinfeld.
It turns out that an appropriate object for which $J$
defines a twisting is a certain weak Hopf algebra that
we describe next.

\subsection{Twisted weak Hopf algebra
$\mathbf{ (\End_k(A)\otimes U)_{J}}$ (cf.\ \cite{Xu1})}
\label{section 3.2}
Observe that the simple algebra $\End_k(A)$ has a natural
structure of a weak Hopf algebra given as follows.

Let $\{ E_{\lambda\mu} \}_{\lambda,\mu\in \mathbb{T}}$ be a basis
of $\End_k(A)$ such that
\begin{equation}
(E_{\lambda\mu}f)(\nu) = \delta_{\mu\nu}f(\lambda), \qquad
f\in A,\, \lambda,\mu,\nu\in \mathbb{T},
\end{equation}
then the comultiplication, counit, and antipode of
$\End_k(A)$ are given by
\begin{equation}
\Delta(E_{\lambda\mu}) = E_{\lambda\mu} \otimes E_{\lambda\mu},
\quad \eps(E_{\lambda\mu}) =1, \quad 
S(E_{\lambda\mu}) = E_{\mu\lambda}.
\end{equation}

\begin{remark}
With the above operations $\End_k(A)$ is a cocommutative weak
Hopf algebra associated to a groupoid with $|\mathbb{T}|$ objects
such that there is a unique morphism $E_{\lambda\mu}$ 
between any two of them, cf.\ Example~\ref{examples}.
\end{remark}

Define the tensor product weak Hopf algebra $H = \End_k(A) \otimes U$.

\begin{proposition}
\label{twist theta}
The elements 
\begin{equation}
\Theta = \sum_{\lambda\mu}\, E_{\lambda\lambda+\mu} \otimes
E_{\lambda\lambda}P_\mu  \quad \mbox{ and } \quad
\bTheta = \sum_{\lambda\mu}\, E_{\lambda+\mu\lambda } \otimes
E_{\lambda\lambda}P_\mu
\label{eqn: Xu's theta}
\end{equation} 
define a twist for  $H$.
\end{proposition}
\begin{proof}
Clearly, we have $\Theta = \Delta(1)\Theta$, $\bTheta = \bTheta\Delta(1)$,
and $\Theta \bTheta =\Delta(1)$. The relations between $\Theta, \bTheta$
and counit are straightforward. We also compute 
\begin{eqnarray*}
(\Delta\otimes\id)(\Theta) (\Theta\otimes 1) &=&
(\id\otimes \Delta)(\Theta) (1\otimes\Theta) \\
&=& \sum_{\lambda\mu\nu }\, E_{\lambda\lambda+\mu+\nu} \otimes
E_{\lambda\lambda+\mu}P_\nu \otimes E_{\lambda\lambda}P_\mu , \\
(\bTheta\otimes 1) (\Delta\otimes\id)(\bTheta) &=&
(1\otimes \bTheta) (\id\otimes \Delta)(\bTheta) \\
&=& \sum_{\lambda\mu\nu }\, E_{\lambda+\mu+\nu\lambda} \otimes
E_{\lambda+\mu \lambda}P_\nu \otimes E_{\lambda\lambda}P_\nu.  
\end{eqnarray*}
One verifies the identities~(\ref{eqn: twist +-})
and (\ref{eqn: twist -+}) of Definition~\ref{twist}
in a similar way.
\end{proof}

Thus, according to Proposition~\ref{twisting},
$H_\Theta = (\End_k(A) \otimes U)_\Theta$ becomes
a weak Hopf algebra. It is non-commutative, non-cocommutative,
and not a Hopf algebra if $|\mathbb{T}|>1$. 
Following \cite{Xu1} we show that it can be further
twisted by means of a dynamical twist $J(\lambda)$ on $U$.

\begin{lemma}
\label{J and Theta}
Let $J(\lambda)\in A \otimes U^{\otimes 2}$
be an invertible, zero weight $U^{\otimes 2}$-valued function on $\mathbb{T}$.
Then the following identities hold in $H^{\otimes 3}$ :
\begin{eqnarray}
(\Delta\otimes \id)(\Theta) (J(\lambda) \otimes 1)
&=& (J(\lambda+h^{(3)}) \otimes 1) (\Delta\otimes \id)(\Theta),  \\
(\id\otimes \Delta)(\Theta)  (1\otimes J(\lambda) )
&=& (1\otimes J(\lambda) ) (\id\otimes \Delta)(\Theta) ,\\
(J^{-1}(\lambda) \otimes 1) (\Delta\otimes \id)(\bTheta)
&=& (\Delta\otimes \id)(\bTheta) (J^{-1}(\lambda+h^{(3)}) \otimes 1),\\
(1\otimes J^{-1}(\lambda) ) (\id\otimes \Delta)(\bTheta)
&=& (\id\otimes \Delta)(\bTheta) (1\otimes J^{-1}(\lambda) ),
\end{eqnarray}
where $J(\lambda) = J\I(\lambda) \otimes J\II(\lambda)$ is embedded
in $H\otimes H$ as 
\begin{equation}
J(\lambda) = \sum_{\lambda} E_{\lambda\lambda}
J\I(\lambda) \otimes E_{\lambda\lambda}J\II(\lambda).
\end{equation}
\end{lemma}
\begin{proof}
We check first two identities, leaving the rest as an exercise
to the reader. Using the formulas for $\Theta$ and $\bTheta$
we compute :
\begin{eqnarray*}
(\Delta\otimes \id)(\Theta) (J(\lambda) \otimes 1)
&=& \sum_{\lambda\mu\nu} 
    E_{\lambda\lambda+\mu}E_{\nu\nu}J\I(\nu) \otimes
    E_{\lambda\lambda+\mu}E_{\nu\nu}J\II(\nu) \otimes
    E_{\lambda\lambda}P_\mu \\
&=& \sum_{\lambda\mu}
    J\I(\lambda+\mu) E_{\lambda\lambda+\mu} \otimes
    J\II(\lambda+\mu) E_{\lambda\lambda+\mu} \otimes
    E_{\lambda\lambda}P_\mu \\
&=& (J(\lambda+h^{(3)}) \otimes 1) (\Delta\otimes \id)(\Theta),  \\
(\id\otimes \Delta)(\Theta)  (1\otimes J(\lambda) )
&=& \sum_{\lambda\mu}
    E_{\lambda\lambda+\mu} \otimes
    E_{\lambda\lambda}{P_\mu}\1 J\I(\lambda) \otimes
    E_{\lambda\lambda}{P_\mu}\2 J\II(\lambda) \\
&=& \sum_{\lambda\mu}
    E_{\lambda\lambda+\mu} \otimes
    J\I(\lambda) E_{\lambda\lambda}{P_\mu}\1  \otimes
    J\II(\lambda) E_{\lambda\lambda}{P_\mu}\2 \\
&=& (1\otimes J(\lambda) ) (\id\otimes \Delta)(\Theta), 
\end{eqnarray*}
where we used the zero weight property of $J(\lambda)$. 
\end{proof}

\begin{proposition}
\label{twist JTheta}
If $J(\lambda)$ is a dynamical twist for $U$ 
then the pair ($F(\lambda), \bF(\lambda)$), where
$$
F(\lambda) = J(\lambda)\Theta \quad \mbox{ and }\quad
\bF(\lambda) = \bTheta J^{-1}(\lambda)
$$
defines a twist for $H= \End_k(A) \otimes U$.
\end{proposition}
\begin{proof}
The twist relations involving counit are obvious,
for the rest we have, using identities from Lemma~\ref{J and Theta} :
\begin{eqnarray*}
(\Delta\otimes\id)(F(\lambda)) (F(\lambda) \otimes 1)
&=& (\Delta\otimes\id)(J(\lambda)\Theta) (J(\lambda)\Theta \otimes 1) \\
&=& (\Delta\otimes\id)(J(\lambda)) J(\lambda+h^{(3)})
    (\Delta\otimes\id)(\Theta) (\Theta \otimes 1) \\
&=& (\id \otimes \Delta)(J(\lambda)) (1\otimes J(\lambda))
     (\id \otimes \Delta)(\Theta) (1\otimes \Theta )\\
&=& (\id \otimes \Delta)(J(\lambda)\Theta ) (1\otimes J(\lambda) \Theta),\\
(\bF(\lambda) \otimes 1) (\Delta\otimes\id)(\bF(\lambda))
&=& (\bTheta J^{-1}(\lambda)\otimes 1) 
     (\Delta\otimes\id)(\bTheta J^{-1}(\lambda)) \\
&=& (\bTheta \otimes 1) (\Delta\otimes\id)(\bTheta ) \\
& & (J^{-1}(\lambda+h^{(3)})\otimes 1)(\Delta\otimes\id)(J^{-1}(\lambda))\\
&=& (1\otimes \bTheta) (\id \otimes \Delta) (\bTheta)
     (1\otimes J^{-1}(\lambda)) (\id \otimes \Delta)(J^{-1}(\lambda))\\
&=& (1\otimes \bF(\lambda)) (\id \otimes \Delta)(\bF(\lambda)), 
\end{eqnarray*}
and one checks other identities similarly.
\end{proof}

Thus, every dynamical twist $J(\lambda)$ for a Hopf algebra $U$
gives rise to a weak Hopf algebra $H_J = H_{J(\lambda)\Theta}$.

\begin{remark}
According to Proposition~\ref{twisting}, the antipode $S_J$ of $H_J$
is given by $S_J(h) = v^{-1}S(h)v$ for all $h\in H_J$, where $S$
is the antipode of $H$ and 
$$
v = \sum_{\lambda \mu} E_{\lambda+\mu \lambda}
    (S(J^{(1)}) J^{(2)})(\lambda) P_\mu.
$$
\end{remark}

\begin{remark}
\label{dynamical gauge transformation}
If $J^x(\lambda)$ is a dynamical twist for $U$ gauge equivalent to
$J(\lambda)$ by means of some $x(\lambda)$ as in 
Remark~\ref{dynamical gauge}, then $x=\sum\, x(\lambda) E_{\lambda \lambda}$
is a gauge transformation of $H$ establishing a gauge equivalence between
the twists ($J(\lambda)\Theta,\,\bTheta J^{-1}(\lambda)$) and
($J^x(\lambda)\Theta,\,\bTheta (J^x)^{-1}(\lambda)$).
\end{remark}

\subsection{Dynamical quantum groups of \cite{EV1}}
\label{section 3.3}
Suppose that $\dim U <\infty$. We will introduce a weak Hopf algebra $D_J$
on a vector space $D = \Map(\mathbb{T}\times\mathbb{T},\, k)\otimes U^*$ 
by adapting the construction of \cite{EV1} to the finite-dimensional case and
show that this weak Hopf algebra is in fact dual to the twisted
weak Hopf algebra $H_J$.

Let $\{ U_i\}$ and $\{ L_i\}$ be dual bases in $U$ and $U^*$,
then the element $L =\sum_i U_i\otimes L_i$ does not depend on the choice
of the bases.
 
Define the coalgebra structure on $D_{J}$ as dual to the algebra
structure of $H$ :
\begin{eqnarray}
(\id \otimes \Delta)L &=& L^{12} L^{13},\\
(\Delta(f))((\lambda^1,\lambda^2),(\mu^1,\mu^2)) 
&=& f(\lambda^1+\mu^1,\lambda^2+\mu^2),\\
(\id\otimes \eps)(L) &=& 1,\\
\eps(f) &=& \sum_\lambda f(\lambda,\lambda), 
\end{eqnarray}
for all $f\in \Map(\mathbb{T}\times\mathbb{T},\,k)$ and 
$\lambda^1,\lambda^2,\mu^1,\mu^2\in \mathbb{T}$.

The algebra structure of $D_{J}$ is given as follows. Observe
that the vector space $U^*$ is bigraded :
$$
U^* =\oplus\, U^*[\alpha^1,\alpha^2],\qquad \mbox { where } \quad
U^*[\alpha^1,\alpha^2] = 
\Hom_k(P_{\alpha^1} U P_{\alpha^2}, k).
$$  
Let us set
\begin{eqnarray}
f(\lambda^1, \lambda^2)g(\lambda^1, \lambda^2) &=& 
g(\lambda^1, \lambda^2)f(\lambda^1, \lambda^2),\\
f(\lambda^1, \lambda^2)L_{\alpha^1 \alpha^2} &=&
L_{\alpha^1 \alpha^2} f(\lambda^1+\alpha^1, \lambda^2+\alpha^2), \\
\label{eqn: Lf=fL}
L^{23}L^{13} &=& : J^{-1}_{12}(\lambda^1)(\Delta \otimes \id)(L) 
J_{12}(\lambda^2):,
\label{eqn: L23 L13}
\end{eqnarray}
where  $f(\lambda^1, \lambda^2), g(\lambda^1, \lambda^2) 
\in \Map(\mathbb{T}\times\mathbb{T},\,k)$, and
$L_{\alpha^1\alpha^2}\in U^*[\alpha^1,\alpha^2]$.
The equation~(\ref{eqn: L23 L13}) is in $U\otimes D_{J}$ and the sign 
$::$ (``normal ordering'') means that the matrix elements of $L$ should be put
on the right of the elements of $J^{-1}(\lambda^1)$, $J(\lambda^2)$.

The unit of $D_{J}$ is defined in an obvious way:
\begin{equation}
1 =(\eps\otimes \id)L.
\end{equation}

Let us consider two  $U$-valued functions on $\mathbb{T}$ :
\begin{equation}
\bK(\lambda) =  m(\id\otimes S)J^{-1}(\lambda), \qquad
K(\lambda) =  (m(S\otimes \id)J)(\lambda-h).
\label{eqn: K K-}
\end{equation}
Note that $K(\lambda)$ and $\bK(\lambda)$ are inverses
of each other and  commute with $A$
since $J(\lambda)$ is of zero weight. Define the antipode
of $D_{J}$ by setting
\begin{eqnarray}
(Sf)(\lambda^1,\lambda^2) &=& f(\lambda^2,\lambda^1), \\
(\id\otimes S)(L) &=& : \bK(\lambda^2)\,(S^{-1}\otimes \id)
(L)\,K(\lambda^1) :, 
\label{eqn: S L}
\end{eqnarray}
for all $f\in \Map(\mathbb{T}\times\mathbb{T},\,k)$ and 
$\lambda^1,\lambda^2\in \mathbb{T}$
and extending $S$ to an algebra anti-homomorphism.

Next, we introduce a basis $\{ \mathbb{I}_{\lambda^1\lambda^2}\}$
of delta-functions on $\mathbb{T}\times\mathbb{T}$, i.e.,
\begin{equation}
\mathbb{I}_{\lambda^1\lambda^2}(\alpha^1,\alpha^2)
= \delta_{\lambda^1 \alpha^1} \delta_{\lambda^2 \alpha^2}.
\end{equation} 
Define a duality form between $D_{J} = \Map(\mathbb{T}\times\mathbb{T},\,k) \otimes U^*
= A\otimes A \otimes U^*$ and $H =\End_k(A) \otimes U$ by
\begin{equation}
\la \mathbb{I}_{\lambda^1\lambda^2} x,\, E_{\alpha^1\alpha^2}u \ra
= \delta_{\lambda^1 \alpha^1} \delta_{\lambda^2 \alpha^2} \la x,\, u\ra
\end{equation}
for all $x\in U^*$ and $u\in U$.

Then in terms of the homogeneous elements
$L_{\alpha^1\alpha^2} \in U^*[\alpha^1,\alpha^2]$ 
the relations (\ref{eqn: L23 L13}) and (\ref{eqn: S L})
defining the multiplication and antipode of $D_{J}$ 
can be rewritten  as
\begin{eqnarray}
\la L_{\alpha^1\alpha^2} L_{\beta^1\beta^2},\, E_{\nu^1 \nu^2}u \ra
&=& \la L_{\beta^1\beta^2} \otimes L_{\alpha^1\alpha^2} ,\,
J^{-1}(\nu^1)\Delta(u) J(\nu^2) \ra 
\label{mult D} \\
\la S(L_{\alpha^1\alpha^2}),\, E_{\nu^1 \nu^2}u \ra
&=& \la L_{\alpha^1\alpha^2},\, \bK(\nu^2) S^{-1}(u) K(\nu^1) \ra,
\label{S D}
\end{eqnarray}
for all $\alpha^1,\alpha^2,\beta^1,\beta^2,\nu^1,\nu^2 \in \mathbb{T}$
and $u\in U$.

\begin{theorem}
\label{duality}
With the above operations $D_{J}$ becomes a weak Hopf algebra
opposite to $H^*_{J}$. 
\end{theorem}
\begin{proof}
We need to show that the structure operations of $D_{J}$
are obtained by transposing the corresponding operations  of $H_{J}$.
For the unit and counit this is obvious.
The product of the elements 
$\mathbb{I}_{\lambda^1\lambda^2}L_{\alpha^1\alpha^2}$
and $\mathbb{I}_{\mu^1\mu^2}L_{\beta^1\beta^2}$ of $D_{J}$, 
where $L_{\alpha^1\alpha^2} \in U^*[\alpha^1,\alpha^2]$, 
can be found by the evaluation against the elements of  $H_{J}$ 
using formula~(\ref{mult D}) :
\begin{eqnarray*}
\lefteqn{ \la (\mathbb{I}_{\lambda^1\lambda^2}L_{\alpha^1\alpha^2}) 
(\mathbb{I}_{\mu^1\mu^2}L_{\beta^1\beta^2}),\, E_{\nu^1\nu^2} u \ra =} \\
&=& \delta_{\lambda^1 \mu^1-\alpha^1} \delta_{\lambda^2 \mu^2-\alpha^2}
    \la \mathbb{I}_{\lambda^1\lambda^2}L_{\alpha^1\alpha^2}L_{\beta^1\beta^2},\,
    E_{\nu^1\nu^2} u \ra  \\
&=& \delta_{\lambda^1 \mu^1-\alpha^1} \delta_{\lambda^1 \nu^1}
    \delta_{\lambda^2 \mu^2-\alpha^2} \delta_{\lambda^2 \nu^2} 
    \la L_{\beta^1\beta^2} \otimes L_{\alpha^1\alpha^2},\, 
    J^{-1}(\nu^1) \Delta(u) J(\nu^2) \ra.
\end{eqnarray*}
On the other hand, we have
\begin{eqnarray*}
\lefteqn{ \la \mathbb{I}_{\mu^1\mu^2}L_{\beta^1\beta^2} \otimes
\mathbb{I}_{\lambda^1\lambda^2}L_{\alpha^1\alpha^2},\, 
\Delta_{J}(E_{\nu^1\nu^2} u) \ra =} \\
&=& \sum_{\eta^1,\eta^2}
   \la \mathbb{I}_{\mu^1\mu^2}L_{\beta^1\beta^2},\, 
   E_{\nu^1+\eta^1\nu^2+\eta^2} \bJ^{(1)}(\nu^1) u\1 J^{(1)}(\nu^2) \ra \\
& & \la \mathbb{I}_{\lambda^1 \lambda^2}L_{\alpha^1 \alpha^2},\, 
  E_{\nu^1 \nu^2} P_{\eta^1} \bJ^{(2)}(\nu^1) u\2 J^{(2)}(\nu^2) P_{\eta^2}\ra\\
&=& \delta_{\lambda^1 \nu^1} \delta_{\lambda^2 \nu^2}
    \la L_{\beta^1 \beta^2},\, \bJ^{(1)}(\nu^1) u\1 J^{(1)}(\nu^2) \ra \\
& & \la L_{\alpha^1 \alpha^2},\,
    P_{\mu^1-\nu^1} \bJ^{(2)}(\nu^1) u\2 J^{(2)}(\nu^2) P_{\mu^2-\nu^2} \ra\\
&=& \delta_{\lambda^1 \nu^1} \delta_{\lambda^2 \nu^2}
    \delta_{\alpha^1 \mu^1-\nu^1} \delta_{\alpha^2  \mu^2-\nu^2} \\
& & \la L_{\beta^1 \beta^2} \otimes L_{\alpha^1 \alpha^2},\,
        J^{-1}(\nu^1) \Delta(u) J(\nu^2) \ra,
\end{eqnarray*}
where $J^{-1} = \bJ^{(1)}\otimes \bJ^{(2)}$.
Comparing the last two relations we conclude that the multiplication
in $D_{J}$ is opposite to the one induced by the comultiplication in
$H_{J}$. 

Finally, the antipode defined by the formula~(\ref{S D})
satisfies 
\begin{eqnarray*}
\lefteqn{ \la S(\mathbb{I}_{\lambda^1\lambda^2}L_{\alpha^1\alpha^2}),\,
E_{\nu^1\nu^2} u \ra =} \\
&=& \la \mathbb{I}_{\lambda^2+\alpha^2 \lambda^1+\alpha^1}S(L_{\alpha^1\alpha^2}),\,
E_{\nu^1\nu^2} u \ra \\
&=& \la \mathbb{I}_{\lambda^2+\alpha^2 \lambda^1+\alpha^1},\,E_{\nu^1\nu^2}\ra
    \la L_{\alpha^1\alpha^2},\,\bK(\nu^2)S^{-1}(u) K(\nu^1) \ra \\
&=& \delta_{\lambda^2+\alpha^2 \nu^1} \delta_{\lambda^1+\alpha^1 \nu^2}
    \la L_{\alpha^1\alpha^2},\,\bK(\nu^2)S^{-1}(u) K(\nu^1) \ra 
\end{eqnarray*}
for all $u\in U$ and $\lambda^1,\lambda^2,\alpha^1,\alpha^2,
\nu^1,\nu^2 \in \mathbb{T}$,
while the transpose of $S_{J}^{-1}$ gives
\begin{eqnarray*}
\lefteqn{ \la \mathbb{I}_{\lambda^1\lambda^2}L_{\alpha^1\alpha^2},\,
S^{-1}(E_{\nu^1\nu^2} u ) \ra = }\\
&=& \sum_{\mu^1\mu^2}\,
    \la \mathbb{I}_{\lambda^1\lambda^2},\, E_{\nu^2+\mu^2  \nu^1-\mu^1} \ra
    \la L_{\alpha^1\alpha^2},\, \bK(\nu^2)P_{-\mu^2} S^{-1}(u) 
       K(\nu^1)P_{\mu^1} \ra \\
&=& \la L_{\alpha^1\alpha^2},\,P_{\nu^2-\lambda^1} \bK(\nu^2)S^{-1}(u)
        K(\nu^1)P_{\nu^1-\lambda^2} \ra \\
&=& \delta_{\alpha^1 \nu^2-\lambda^1} \delta_{\alpha^2 \nu^1-\lambda^2}
    \la L_{\alpha^1\alpha^2},\,\bK(\nu^2) S^{-1}(u) K(\nu^1) \ra,
\end{eqnarray*}
which completes the proof.
\end{proof}

\subsection{The case of quasitriangular $\mathbf{U}$.}
\label{section 3.4}
When $U$ is quasitriangular with the universal $R$-matrix  $\mR$, 
the twisted weak Hopf algebra $H_{J}$ is also quasitriangular 
by means of the matrices
\begin{equation}
\mR(\lambda) = \bTheta J_{21}^{-1}(\lambda) \mR J(\lambda) \Theta
\quad \mbox{ and } \quad
\bmR(\lambda) = \bTheta J^{-1}(\lambda) \bmR J_{21}(\lambda) \Theta.
\end{equation}
More explicitly,
\begin{equation}
\label{eqn: twisted R-matrix}
\mR(\lambda) = \sum_{\lambda\mu\nu}\, E_{\lambda \lambda+\nu}P_\mu
\mR^{J(1)}(\lambda) \otimes E_{\lambda+\mu \lambda} \mR^{J(2)}(\lambda)P_\nu,
\end{equation}
where $\mR^{J}(\lambda) =  J_{21}^{-1}(\lambda) \mR J(\lambda)$. 
By~(\ref{eqn: H* to H}), $\mR(\lambda)$ establishes
a weak Hopf algebra homomorphism  $H^*_{J} \to H^{op}_{J}$,
i.e., a homomorphism $\rho: D_{J} \to H_{J}$. 

\begin{remark}
\label{NZ}
Note that $\dim(P_\mu U) = \dim( U P_\nu) = \dim U/ |\mathbb{T}|$
for all $\mu,\nu\in \mathbb{T}$,
since any finite dimensional Hopf algebra is free
over its Hopf subalgebra  \cite{M}. In particular, 
$\dim U/ |\mathbb{T}|$ is an integer.
\end{remark}

\begin{proposition}
\label{non-degeneracy criterion}
$\rho$ is an isomorphism if and only if the element
\begin{equation}
\mR_{\mu\nu}^{J}(\lambda) =
P_\mu \mR^{J(1)}(\lambda) \otimes \mR^{J(2)}(\lambda) P_\nu
\in U\otimes U 
\label{eqn: R mu nu}
\end{equation}
has the maximal possible rank $(= \dim U/ |\mathbb{T}|)$ for all fixed
$\lambda,\mu,\nu\in \mathbb{T}$.
\end{proposition}
\begin{proof}
Since $H$ and $D$ are finite dimensional, $\rho$ is an isomorphism
if and only if its image coincides with $H$, i.e.,
\begin{eqnarray*}
\Image(\rho) 
&=&  \span \{ \sum_\nu E_{\lambda \lambda+\nu}P_\mu \mR^{J(1)}(\lambda) 
     \phi( \mR^{J(2)}(\lambda) P_\nu ) 
     \mid \lambda,\mu\in \mathbb{T}, \phi\in U^* \}  \\
&=& \bigoplus_{\lambda\mu\nu} E_{\lambda \lambda+\nu}
    \span \{ P_\mu \mR^{J(1)}(\lambda)\, \phi(\mR^{J(2)}(\lambda) P_\nu) 
    \mid \phi\in U^* \} = U,
\end{eqnarray*}
which happens precisely when
\begin{equation}
\span \{ P_\mu \mR^{J(1)}(\lambda)\, \phi(\mR^{J(2)}(\lambda) P_\nu)
\mid \phi\in U^* \} = P_\mu U,
\end{equation}
for all $\lambda,\mu,\nu\in \mathbb{T}$. Since a Hopf algebra $U$ is a free 
$\Map(\mathbb{T},\,k)$-module, we see that the rank  of
$\mR_{\mu\nu}^{J}(\lambda)$ has to be equal to
$\dim(P_\mu U)$ and, therefore, to $\dim{U}/|\mathbb{T}|$, by
Remark~\ref{NZ}. 
\end{proof}

\end{section}  
  

\begin{section}
{Dynamical twists for $U_q(\g)$ at roots of $1$} 
\label{4}

\subsection{Construction of $\mathbf{ J(\lambda)}$}
\label{section 4.1}
Suppose that $\g$ is a simple Lie algebra of type
$A$, $D$ or $E$ and $q$ is a primitive $\ell$th root of unity in $k$,
where $\ell\geq 3$ is odd and coprime with the determinant
of the Cartan matrix $(a_{ij})_{ij=1,\dots,m}$ of $\g$.

Let $U=U_q(\g)$ be the corresponding quantum group
which is a finite dimensional Hopf algebra with generators
$E_i,F_i, K_i$, where $i=1,\dots,m$ and the following relations \cite{L} :
\begin{eqnarray*}
& & K_i^\ell =1,\quad E_i^\ell =0,\quad F_i^\ell =0, \\
& & K_iK_j = K_jK_i,\quad K_iE_j = q^{a_{ij}}E_jK_i,\quad
K_iF_j = q^{-a_{ij}}F_j K_i, \\
& & E_iF_j - F_jE_i =\delta_{ij} \dfrac{K_i-K_i^{-1}}{q-q^{-1}},\\
& & E_iE_j = E_jE_i,\quad F_iF_j = F_jF_i \quad \mbox{ if } \quad a_{ij}=0,\\
& & E_i^2E_j -(q+q^{-1})E_iE_jE_i + E_jE_i^2 = 0  
    \quad \mbox{ if }\quad a_{ij}=-1,\\ 
& & F_i^2F_j -(q+q^{-1})F_iF_jF_i + F_jF_i^2 = 0 
   \qquad \mbox{ if } \quad a_{ij}=-1,
\end{eqnarray*}
with the comultiplication, counit, and antipode given by
\begin{eqnarray*}
\Delta(K_i) &=& K_i\otimes K_i, \quad
    \Delta(E_i) = E_i \otimes 1 + K_i\otimes E_i, \quad
    \Delta(F_i) = F_i \otimes K_i^{-1} + 1\otimes F_i, \\ 
S(K_i) &=& K_i^{-1}, \quad S(E_i) = -K_i^{-1}E_i, \quad
    S(F_i) = - F_i K_i, \\
\eps(K_i) &=& 1,  \quad \mbox{ and } \quad \eps(E_i) =\eps(F_i) =0. 
\end{eqnarray*}
Denote by $(\cdot,\cdot)$ the bilinear form on $\mathbb{Z}^m$
defined by the Cartan matrix $(a_{ij})$. 

Let $\mathbb{T} \cong (\mathbb{Z}/\ell\mathbb{Z})^m$ be the Abelian
group generated by $K_i,\,i=1,\dots,m$. For any $m$-tuple of integers 
$\gamma =(\gamma_1,\dots,\gamma_m)$
we will write $K_\gamma = K_1^{\gamma_1}\dots K_m^{\gamma_m}\in \mathbb{T}$. 
Denote by $I$ the set of all integer $m$-tuples
$\gamma =(\gamma_1,\dots,\gamma_m)$ with $0\leq \gamma_i< \ell,\,
i=1,\dots,m$.  

Let $\Delta^+$ be the set
of positive roots of $\g$. For each $\alpha\in \Delta^+$
define $E_\alpha,F_\alpha\in U_q(\g)$ inductively
by setting $E_{\alpha_i}=E_i,\, F_{\alpha_i}=F_i$ for all
simple roots $\alpha_i,\,i=1,\dots,m$ and
\begin{equation*}
E_{\alpha+\alpha'} =q^{-1}E_{\alpha'}E_\alpha - E_\alpha E_{\alpha'}
\quad \mbox{and} \quad
F_{\alpha+\alpha'} =qF_\alpha F_{\alpha'} - F_{\alpha'} F_\alpha,
\qquad \alpha,\alpha'\in \Delta^+.
\end{equation*} 
Let $\beta_1,\dots,\beta_N$ be the normal ordering of
$\Delta^+$ and for every  $N$-tuple of non-negative integers
$a=(a_1,\dots, a_N)$ introduce the monomials
\begin{equation*}
E_a = E_{\beta_1}^{a_1}\dots E_{\beta_N}^{a_N}
\quad \mbox{ and } \quad
F_a = F_{\beta_1}^{a_1}\dots F_{\beta_N}^{a_N}.
\end{equation*}
Then the universal $R$-matrix of $U_q(\g)$ is given by \cite{R}, \cite{T} :
\begin{equation}
\label{eqn: R-matrix}
\mR = \frac{1}{\ell^m} \prod_{s=1}^{N} \left( \sum_{n=0}^{\ell-1}
q^{-\frac{n(n+1)}{2} } \frac{(1-q^2)^n}{[n]_q!} E_{\beta_s}^n
\otimes F_{\beta_s}^n
\right) \left( \sum_{\beta,\gamma\in I} q^{(\beta,\gamma)}
K_\beta \otimes  K_\gamma \right),
\end{equation}
where $[n]_q! = [1]_q[2]_q\cdots [n]_q$, 
$[n]_q= \frac{q^n-q^{-n}}{q-q^{-1}}$, and  
\begin{equation}
\Omega =  \frac{1}{\ell^m} \sum_{\beta,\gamma\in I} q^{(\beta,\gamma)}
K_\beta \otimes  K_\gamma,
\end{equation}
is the ``Cartan part'' of $\mR$. 

Note that the idempotents
\begin{equation*}
P_\beta = \frac{1}{|\mathbb{T}|} \sum_{\lambda\in I}\, 
q^{(\beta,\lambda)}K_\lambda
\end{equation*}
generate a commutative and cocommutative Hopf subalgebra 
$A =\Map(\mathbb{T},k)$ of $U$.

Observe that $U$ is $\mathbb{Z}^m$-graded in such a way that
a monomial $X$ in $E_i,F_i, K_i$
belongs to $U[\beta-\beta']$ where $\beta =(\beta_1,\dots,\beta_m)$
and  $\beta' =(\beta_1',\dots,\beta_m')$ are such that
$E_i$ appears exactly $\beta_i$ times and
$F_i$ appears exactly $\beta_i'$ times in $X$ for each $i$.

This induces a $\mathbb{Z}$-grading of the algebra $U$ with
\begin{equation}
\deg(E_i) =1, \quad \deg(F_i) = -1, \quad \deg(K_i) = 0,
\qquad i=1,\dots,m, 
\end{equation}
and $\deg(XY)= \deg(X)+\deg(Y)$ for all $X$ and $Y$. Of course,
there are only finitely many non-zero components of $U$
since it is finite dimensional.

Let $U_+$ be the subalgebra of $U$ generated by
the elements $E_i,\,K_i,\,i=1,\dots,m$, $U_-$ be
the subalgebra  generated by $F_i,\,K_i,\,i=1,\dots,m$, and
$I_\pm$ be the kernels of the projections from $U_\pm$
to the elements of zero degree.

For arbitrary non-zero constants $\Lambda_1,\dots, \Lambda_m$
define a Hopf algebra automorphism $\Lambda$ of $U$ by
setting 
\begin{equation*}
\Lambda(E_i) = \Lambda_i E_i,\quad 
\Lambda(F_i) = \Lambda_i^{-1} F_i,\quad  \mbox{ and } \quad
\Lambda(K_i) = K_i \quad\mbox{ for all } i=1,\dots,m.
\end{equation*}
If for $\beta =(\beta_1,\dots,\beta_m) \in \mathbb{Z}^m$ 
we write $\Lambda_\beta = \Lambda_1^{\beta_1}\dots
\Lambda_m^{\beta_m}$ then
\begin{equation}
\Lambda|_{U[\beta]} = \Lambda_\beta\, \id_{U[\beta]}.
\end{equation}

\begin{definition}
\label{generic Lambda}
We will say that $\Lambda=(\Lambda_1,\dots,\Lambda_m)$ is {\em generic} 
if the spectrum of $\Lambda$ does not contain $\ell$th roots of unity.
\end{definition}

For every $K_\lambda\in \mathbb{T}$ we introduce the following linear
operator on $U\otimes U$:
\begin{equation}
\label{AL2}
A_2^L(\lambda) X =(\Ad K_\lambda \circ \Lambda \otimes \id) (\mR X\Omega^{-1}).
\end{equation}

\begin{proposition}
\label{equation for J}
For every generic $\Lambda$ there exists a unique element
$J(\lambda) \in 1 + I_+ \otimes I_-$ that satisfies the
following ABRR relation \cite{ABRR}, \cite{ES1}, \cite{ESS} :
\begin{equation}
\label{ABRR}
A_2^L(\lambda) J(\lambda) = J(\lambda).
\end{equation}
\end{proposition}
\begin{proof}
Let us write $X =\sum_{j\geq 1}\, X^j$, where $X^j$ is
the sum of all terms having the $\mathbb{Z}$-degree $j$ 
in the first component. Using the structure of the $R$-matrix
of $U$, we can write (\ref{ABRR}) as a finite system of linear
equations labeled by degree $j\geq 1$:
\begin{equation}
\label{recursive system}
X^j = (\Ad K_\lambda \circ \Lambda \otimes \id) (\Omega X^j \Omega^{-1})
+ \cdots,
\end{equation}
where $\cdots$ stands for the terms involving $X^i$ for $i<j$.
Thus, the equation (\ref{ABRR}) can be solved recursively,
starting with $X^0 =1$, and the solution is unique provided
that the operator 
\begin{equation}
\label{operator}
\id - (\Ad K_\lambda \circ \Lambda \otimes \id)\circ \Ad \Omega
\end{equation}
is invertible in $\End_k(U\otimes U)$. Let us show that this operator
is diagonalizable and compute its eigenvalues.

For all $X_\alpha\in U[\alpha]$, $X_{\alpha'}\in U[\alpha']$,
and $K_\beta\in \mathbb{T}$ we have
$$
(\Ad K_\beta)X_\alpha = q^{(\alpha,\beta)} X_\alpha ,
\qquad
(\Ad K_\beta)X_{\alpha'} = q^{(\alpha',\beta)} X_{\alpha'},
$$
and therefore
$$
\Ad \Omega(X_\alpha \otimes X_{\alpha'})
= q^{(\alpha,\alpha')} (K_{-\alpha'}\otimes K_{-\alpha})
(X_\alpha \otimes X_{\alpha'}),
$$
whence the eigenvalues of $\Ad \Omega$ in $U[\alpha]\otimes U[\alpha']$
are the numbers
$$
d_{\chi\chi'} = q^{(\alpha,\alpha')}  \chi(K_{-\alpha}) \chi'(K_{-\alpha'}),
$$
where $\chi$ and $\chi'$ are characters of $\mathbb{T}$.
In particular, each $d_{\chi\chi'}$ is an $\ell$th root of unity.

Clearly, the eigenvalue of the operator 
$(\Ad K_\lambda \circ \Lambda \otimes \id)$ in 
$U[\alpha]\otimes U[\alpha']$ is $\Lambda_\alpha q^{(\alpha,\lambda)}$.
Putting these numbers together, we conclude that (\ref{operator})
is invertible in $U\otimes U$ if and only if
$$
1- \Lambda_\alpha q^{(\alpha,\lambda+\alpha')} d_{\chi\chi'} \neq 0,
$$
for all $\lambda,\alpha,\alpha',\chi,\chi'$ which is
the case for generic $\Lambda$.
\end{proof}

\begin{remark}
\label{J commutes with Ad K}
If $J(\lambda)\in 1+ I_+\otimes I_-$ is a solution of (\ref{ABRR}),
then it has zero weight,
by the uniqueness result of 
Proposition~\ref{equation for J}, since $(\Ad \Delta(K_\beta)) J(\lambda)$
is also a solution of (\ref{ABRR}) for all $K_\beta \in \mathbb{T}$.
\end{remark}

Our goal is to show that the above element $J(\lambda)$
gives rise to a dynamical twist for $U_q(\g)$.

Similarly to (\ref{AL2}) define the operator
\begin{equation}
\label{AR2}
A^2_R(\lambda) X 
= (\id \otimes \Ad K_{-\lambda}\circ \Lambda^{-1})(\mR X\Omega^{-1})
\end{equation}
for all $X\in U\otimes U$.

\begin{lemma}
\label{AL2 and AR2} 
$A_L^2(\lambda)$ and $A^2_R(\lambda)$ commute.
\end{lemma}
\begin{proof}
For all $X\in U\otimes U$ we have
\begin{eqnarray*}
\lefteqn{ (\Ad K_\lambda\circ \Lambda \otimes \id)\mR\,
    (\Ad K_\lambda\circ \Lambda \otimes \Ad K_{-\lambda}\circ \Lambda^{-1})
    (\mR X\Omega^{-1})\,\Omega^{-1} = }\\
&=& (\id \otimes \Ad K_{-\lambda}\circ \Lambda^{-1})\mR\,    
    (\Ad K_\lambda\circ \Lambda \otimes \Ad K_{-\lambda}\circ \Lambda^{-1})
    (\mR X\Omega^{-1})\,\Omega^{-1}, 
\end{eqnarray*}
since $(\Lambda \otimes \id)\mR = (\id \otimes \Lambda^{-1})\mR$,
whence $A_L^2(\lambda)\circ A_R^2(\lambda) =   
A_R^2(\lambda) \circ A_L^2(\lambda)$.
\end{proof}

\begin{corollary}
\label{J solves the system}
$J(\lambda)$ is the unique solution of the system
\begin{equation}
A_L^2(\lambda) X = X  \qquad \mbox{ and } \qquad
A_R^2(\lambda) X = X.
\end{equation}
with $X\in 1 + I_+\otimes I_-$.
\end{corollary}
\begin{proof}
We have 
$$
A_L^2(\lambda) \circ A_R^2(\lambda) J(\lambda) = 
A_R^2(\lambda) \circ  A_L^2(\lambda) J(\lambda)  =   
A_R^2(\lambda) J(\lambda),
$$ 
hence $ A_R^2(\lambda) J(\lambda) = J(\lambda)$ by the uniqueness 
of the solution  of (\ref{ABRR}). 
\end{proof}

Following \cite{ESS}, consider the $3$-component operators :
\begin{eqnarray}
A_L^3(\lambda) X 
&=& (\Ad K_\lambda\circ \Lambda \otimes \id \otimes \id)
(\mR_{13}\mR_{12} X \Omega_{12}^{-1} \Omega_{13}^{-1}) 
\label{AL3} \\
A_R^3(\lambda) X 
&=& (\id \otimes \id \otimes \Ad K_{-\lambda}\circ \Lambda^{-1})
(\mR_{13}\mR_{23} X \Omega_{13}^{-1} \Omega_{23}^{-1})
\label{AR3}
\end{eqnarray}

\begin{lemma}
\label{AL3 and AR3 commute}
The operators $A_L^3(\lambda)$ and $A_R^3(\lambda)$ commute.
\end{lemma}
\begin{proof}
This statement amounts to showing that
\begin{eqnarray*}
\lefteqn{ (\Ad K_\lambda\circ \Lambda \otimes \id \otimes \id) \mR_{13}\mR_{12}\,
(\Ad K_\lambda\circ \Lambda \otimes \id \otimes \Ad K_{-\lambda}\circ \Lambda^{-1}) 
\mR_{13}\mR_{23}=}\\
&=& (\id \otimes \id \otimes \Ad K_{-\lambda}\circ \Lambda^{-1})\mR_{13}\mR_{23}\,
(\Ad K_\lambda\circ \Lambda \otimes \id \otimes \Ad K_{-\lambda}\circ \Lambda^{-1}) 
\mR_{13}\mR_{12}.
\end{eqnarray*}
If we denote $\hat{\mR} =(\Ad K_\lambda\circ \Lambda \otimes \id )\mR
= (\id \otimes \Ad K_{-\lambda}\circ \Lambda^{-1})\mR$ and
$$
\tilde{\mR}= (\Ad K_\lambda\circ \Lambda \otimes
\Ad K_{-\lambda}\circ \Lambda^{-1})\mR,
$$  
then the above equality translates to
\begin{equation*}
\hat{\mR}_{13}\hat{\mR}_{12}\tilde{\mR}_{13}\hat{\mR}_{23}
= \hat{\mR}_{13} \hat{\mR}_{23} \tilde{\mR}_{13} \hat{\mR}_{12},
\end{equation*}
which follows from the quantum Yang-Baxter equation
after cancelling the first factor.
\end{proof}

\begin{lemma}
\label{AL3 X = AR3 X =X}
If there exists a solution $X$ of the system
\begin{equation}
A_L^3(\lambda) X = A_R^3(\lambda) X =X
\end{equation}
with $X \in I_+\otimes U\otimes U + U\otimes U\otimes I_-$,
then it is unique.
\end{lemma}
\begin{proof}
It is enough to show that such a solution $X$ is unique
for the equation $A_L^3 A_R^3 X =X$. Let us write
$$
X = 1 +\sum_{k,l\geq 0; k+l>0}\, X^{k,l}
$$
where $X^{k,l}$ is the sum of all elements having $\mathbb{Z}$-degree
$k$ in the first component and $-l$ in the third one.
Then the equation $A_L^3(\lambda) A_R^3(\lambda) X =X$ 
transforms to the system
\begin{eqnarray*}
X^{k,l} 
&=& (\Ad K_\lambda\circ \Lambda \otimes \id 
    \otimes \Ad K_{-\lambda}\circ \Lambda^{-1})
    (W X^{k,l} W^{-1}) \\*
&  & + \mbox{ terms depending on } X^{k',l'} \mbox{ with }
     k'+l' < k+l,
\end{eqnarray*}
for all $k\geq 0$ and $l\geq 0$ such that $k+l>0$, where
$W = \Omega_{12} \Omega_{23} (\Omega_{13})^2$.

As in Proposition~\ref{equation for J} one can check
that the operator
\begin{equation}
\id - (\Ad K_\lambda\circ \Lambda \otimes  \id \otimes 
\Ad K_{-\lambda}\circ \Lambda^{-1}) \circ \Ad W
\end{equation}
is invertible for generic $\Lambda$, therefore the above system
can be solved recursively and the solution is unique.
\end{proof}

\begin{theorem}
\label{J is a shifted twist}
$J(\lambda)$ satisfies the equations
\begin{eqnarray}
(\Delta\otimes\id)J(\lambda) ( J(\lambda+h^{(3)})\otimes 1)
&=& (\id\otimes\Delta)J(\lambda) (1 \otimes J(\lambda-h^{(1)})), 
\label{eqn: shifted twist} \\
(\eps\otimes\id)J(\lambda) &=& (\id\otimes\eps)J(\lambda) =1.
\end{eqnarray}
\end{theorem}
\begin{proof}
Let us denote the left-hand side of (\ref{eqn: shifted twist})
by $Y_L$ and the right-hand side by $Y_R$. We show that both $Y_L$ and $Y_R$
are solutions of the system $A_L^3(\lambda) X = A_R^3(\lambda) X = X$, 
then the result will follow from Lemma~\ref{AL3 X = AR3 X =X}. We have :
\begin{eqnarray*}
A_R^3 Y_L
&=& (\id \otimes \id \otimes \Ad K_{-\lambda}\circ \Lambda^{-1})
    (\mR_{13}\mR_{23} Y_L \Omega_{13}^{-1} \Omega_{23}^{-1}) \\
&=& (\id \otimes \id \otimes \Ad K_{-\lambda}\circ \Lambda^{-1})
    (\Delta\otimes \id)(\mR J(\lambda))  J(\lambda+h^{(3)})
    (\Delta\otimes \id)\Omega^{-1}\\
&=& (\Delta\otimes \id)(\id \otimes \Ad K_{-\lambda}\circ \Lambda^{-1})
    (\mR J(\lambda)\Omega^{-1}) (J(\lambda+h^{(3)}) \otimes 1)\\
&=& (\Delta\otimes \id)J(\lambda) (J(\lambda+h^{(3)}) \otimes 1) =Y_L, \\
A_L^3 Y_R
&=& (\Ad K_\lambda\circ \Lambda \otimes \id \otimes \id)
    (\mR_{13}\mR_{12} Y_R \Omega_{12}^{-1} \Omega_{13}^{-1}) \\
&=& (\Ad K_\lambda\circ \Lambda \otimes \id \otimes \id)
     (\id\otimes\Delta)(\mR J(\lambda)) (1 \otimes J(\lambda-h^{(1)}))
     (\id\otimes\Delta)\Omega^{-1} \\
&=& (\id\otimes\Delta) (\Ad K_\lambda\circ \Lambda \otimes \id )
    (\mR J(\lambda) \Omega^{-1}) (1 \otimes J(\lambda-h^{(1)})) \\
&=& (\id\otimes\Delta) J(\lambda) (1 \otimes J(\lambda-h^{(1)})) =Y_R,
\end{eqnarray*}
where we used that $A_L^2 J = A_R^2 J =J$, properties of the
$R$-matrix, and that $J$ commutes with $\Ad K_\lambda$.

To establish that $A_L^3 Y_L =Y_L$ note that both $A_L^3 Y_L$ and $Y_L$
are solutions of the equation $A_R^3 X =X$ and therefore are determined
uniquely by their parts of zero degree in the third component.
Thus it suffices to compare these parts :
\begin{eqnarray*}
(A_L^3 Y_L)^0
&=& (\Ad K_\lambda\circ \Lambda \otimes \id \otimes \id)\circ \Ad \Omega_{13}
    (\mR_{12} J_{12}(\lambda+h^{(3)}) \Omega_{12}^{-1}) \\
&=& \sum_{\beta} (\Ad K_{\lambda+\beta} \circ \Lambda \otimes \id )
    (\mR J(\lambda+h^{(3)})\Omega^{-1}) \otimes  P_\beta \\
&=& (\Ad K_{\lambda+h^{(3)}} \circ \Lambda \otimes \id )
    (\mR J(\lambda+h^{(3)})\Omega^{-1}) = J(\lambda+h^{(3)}) = Y_L^0,\\
(A_R^3 Y_R)^0
&=& (\id \otimes \id \otimes \Ad K_{-\lambda}\circ \Lambda^{-1})
    \circ\Ad \Omega_{13} (\mR_{23} J_{23}(\lambda-h^{(1)})\Omega_{23}^{-1}) \\
&=& \sum_{\beta} 
    P_\beta \otimes  (\id \otimes \Ad K_{-\lambda+\beta}\circ \Lambda^{-1})
    (\mR J(\lambda-h^{(1)}) \Omega^{-1}) \\
&=& (\id \otimes \Ad K_{-(\lambda-h^{(1)})} \circ \Lambda^{-1})
    (\mR J(\lambda-h^{(1)}) \Omega^{-1})  = J(\lambda-h^{(1)}) = Y_R^0.
\end{eqnarray*}
The relations between $J(\lambda)$ and $\eps$ are obvious. 
\end{proof}

\begin{proposition}
\label{shift undone}
The element  
\begin{equation}
\tJ(\lambda)= J(2\lambda+h^{(1)}+h^{(2)})
\end{equation} 
is a dynamical twist for $U_q(\g)$
in the sense of Definition~\ref{dynamical twist}. 
\end{proposition}
\begin{proof}
We directly  compute :
\begin{eqnarray*}
\lefteqn{ 
(\Delta\otimes\id)\tJ(\lambda) (\tJ(\lambda+h^{(3)})\otimes 1)=}\\
&=& (\Delta\otimes\id)J(2\lambda+h^{(1)}+h^{(2)}+h^{(3)})
    (J(2\lambda+h^{(1)}+h^{(2)}+2h^{(3)}) \otimes 1) \\
&=& (\id\otimes \Delta)J(2\lambda+h^{(1)}+h^{(2)}+h^{(3)})
    (1\otimes J(2\lambda+h^{(2)}+h^{(3)})) \\
&=& (\id\otimes \Delta)\tJ(\lambda) (1\otimes \tJ(\lambda)).  
\end{eqnarray*} 
The identities $(\eps\otimes\id) \tJ(\lambda) =1$ and 
$(\id \otimes \eps)\tJ(\lambda) =1$ are clear. 
\end{proof}

\begin{example}
\label{sl 2}
Let us give an explicit expression for the twists 
$J(\lambda)$ and $\tJ(\lambda)$ in the case 
$\g=sl(2)$. In this case, $U_q(\g)$ is generated by $E,F,K$ with the standard 
relations. The element analogous to 
$J(\lambda)$ for generic $q$ was computed 
already in \cite{Ba} (see also \cite{BBB}). 
If we switch to our conventions, this element will take the form
$$
J(\lambda)=\sum_{n=0}^\infty q^{-n(n+1)/2}\frac{(1-q^2)^n}{[n]_q!}
(E^n\otimes F^n)
\prod_{\nu=1}^n\frac{ \Lambda q^{2\lambda}}{1-\Lambda q^{2\lambda+2\nu}
(K\otimes K^{-1})}.
$$
It is obvious that the formula for 
$q$ being a primitive $\ell$-th root of unity is simply obtained 
by truncating this formula:
$$
J(\lambda)=\sum_{n=0}^{\ell-1} q^{-n(n+1)/2}\frac{(1-q^2)^n}{[n]_q!}
(E^n\otimes F^n)
\prod_{\nu=1}^n\frac{ \Lambda q^{2\lambda}}{1-\Lambda q^{2\lambda+2\nu}
(K\otimes K^{-1})}.
$$
Therefore, 
$$
\tJ(\lambda)=\sum_{n=0}^{\ell-1} q^{-n(n+1)/2}\frac{(1-q^2)^n}{[n]_q!}
(E^n\otimes F^n)
\prod_{\nu=1}^n\frac{\Lambda q^{4\lambda}K\otimes K}{1-\Lambda q^{4\lambda+2\nu}
(K^2\otimes 1)}.
$$ 
Note that the term of this sum corresponding to $n=1$ coincides with the one 
computed in Section~\ref{section 4.3}.
\end{example}

\subsection{Dynamical twists arising from generalized
Belavin-Drinfeld triples (cf. \cite{ESS})}
\label{section 4.2}

\begin{definition}
\label{BD triple}
A generalized Belavin-Drinfeld triple for a simple
Lie algebra $\g$ consists of subsets
$\Gamma_1,\Gamma_2$ of the set 
$\Gamma=(\alpha_1,\dots,\alpha_m)$ of simple roots 
of $\g$ together with an inner product preserving
bijection $T: \Gamma_1 \to \Gamma_2$.
\end{definition}

We say that $T$ is {\em nilpotent} if for any $i=1,\dots,m$
there exists a positive integer $d_i$ such that $T^{d_i}(\alpha_i) 
\not\in \Gamma_1$.
For non-nilpotent $T$ we define an {\em order} of $T$ (denoted by $n(T)$) 
to be the least common multiple of the lengths of orbits of $T$.

Let $Q_1,\, Q_2$, and $Q$ be the free Abelian groups generated
by the root sets $\Gamma_1$, $\Gamma_2$, and $\Gamma$ 
respectively and
\begin{equation}
L = \{ \lambda\in Q \mid (\lambda,\,\alpha) = (\lambda,\,T\alpha) \quad 
\forall \alpha \in \Gamma_1\}.
\end{equation}
Then $T$ extends to the isomorphism between $Q_1$ and $Q_2$ and one can
check as in Lemma~3.1 in \cite{ESS}
that  $Q_1\cap L = \{ \lambda\in Q_1 \mid T\lambda = \lambda \}$
and $Q_1+L$ and $L^\perp+L$ are finite index subgroups of $Q$. Let
\begin{equation}
n_1 = [Q: (Q_1+L)] \quad \mbox{ and } \quad n_2 = [Q: (L^\perp+L)] 
\end{equation}
We assume that $\ell$ is coprime with both $n_1$ and $n_2$ and
introduce a homomorphism $Q_1+L \to Q$ also denoted by $T$, letting
\begin{equation}
T|_{Q_1} = T, \qquad T|_L = \id.
\end{equation}
Since $T=\id$ on $Q_1\cap L$ it follows that $T$ is well-defined.

Factorization by $\ell Q$ yields an automorphism of $\mathbb{T} =Q/\ell Q$ :
\begin{equation}
\mathcal{T} \quad : \quad (Q_1+L)/\ell Q = Q/\ell Q \to Q/\ell Q.
\end{equation}
which preserves the inner product $(\cdot,\cdot) : \mathbb{T}\times \mathbb{T} 
\to  \mathbb{Z}/\ell \mathbb{Z}$ and extends to the 
algebra homomorphisms $\mathcal{T}_\pm: U_\pm \to U_\pm$ defined by
\begin{eqnarray}
\mathcal{T}_+(K_\lambda) &=& K_{T\lambda}, \quad
\mathcal{T}_-(K_\lambda) = K_{T^{-1}\lambda}, \\
\mathcal{T}_+(E_i) &=& E_{T(i)} \quad \mbox{ if } \alpha_i\in \Gamma_1
\mbox{ and } \mathcal{T}_+(E_i)=0 \mbox{ otherwise }, \\
\mathcal{T}_-(F_i) &=& F_{T^{-1}(i)} \mbox{ if } \alpha_i\in \Gamma_2
\mbox{ and } \mathcal{T}_-(F_i)=0 \mbox{ otherwise },  
\end{eqnarray}
for all $\lambda \in \mathbb{T}$ and $i=1,\dots, m$, where
$T(i)$ denotes the number such that $T(\alpha_i) = \alpha_{T(i)}$
and $T^{-1}(i)$ the number such that  $T^{-1}(\alpha_i) = \alpha_{T^{-1}(i)}$.

Let $\mathbb{T}_L$ and $\mathbb{T}_L^\perp$
be the images of $L$ and $L^\perp$ in $\mathbb{T}$ and
\begin{equation}
\Omega_L = \sum_{K_\beta,K_\gamma\in \mathbb{T}_L}\,
q^{(\beta,\gamma)} K_\beta\otimes K_\gamma,
\qquad
\Omega_{L^\perp} = \sum_{K_\beta,K_\gamma\in \mathbb{T}_L^\perp}\,
q^{(\beta,\gamma)} K_\beta\otimes K_\gamma.
\end{equation}
Then $\mathbb{T} = \mathbb{T}_L \oplus \mathbb{T}_L^\perp$ and
$\Omega =  \Omega_L \Omega_{L^\perp}$.
We define a modification of the operator $A_L^2(\lambda)$
introduced in (\ref{AL2}):
\begin{equation}
A_L^2(\lambda)X = (\mathcal{T}_+\circ \Ad K_\lambda\circ \Lambda \otimes \id)
(\mR X \Omega_L^{-1}), \qquad \lambda \in \mathbb{T}_L.
\end{equation}
One can show, using the same argument as in 
Proposition~\ref{J commutes with Ad K}, 
that if  $\Lambda_i = \Lambda_{T(i)}$ for all $\alpha_i\in \Gamma_1$ and
the spectrum of $\Lambda$ (in the case $T$ is not nilpotent)
does not contain roots of unity of order $n(T)\ell$,
then there exists a unique element 
$J_T(\lambda) \in Z + I_+ \otimes I_-$,
where $Z=((\id-T)^{-1}T \otimes \id) \Omega_{L^\perp}$,
satisfying the modified ABRR relation (cf.\ \cite{ESS}):
\begin{equation}
A_L^2(\lambda) J_T(\lambda) = J_T(\lambda)
\end{equation}
and commuting with $\Ad \Delta(K_\lambda)$ for all $\lambda\in \mathbb{T}_L$.
Next, if we define
\begin{equation}
A_R^2(\lambda)X = (\id \otimes \mathcal{T}_-\circ \Ad K_\lambda\circ \Lambda)
(\mR X \Omega_L^{-1}), \qquad \lambda \in \mathbb{T}_L,
\end{equation}
then Lemma~\ref{AL2 and AR2} and Corollary~\ref{J solves the system}
are still valid because of the identity
$$
(\mathcal{T}_+ \otimes \id)\Omega_L =
(\id \otimes \mathcal{T}_-)\Omega_L 
$$
that follows from the inner product preserving
property of $T$.

Modifying definitions of $A_L^3(\lambda)$ and $A_R^3(\lambda)$
as in \cite{ESS}, it is possible to repeat the proofs of 
Lemmas~\ref{AL3 and AR3 commute}, \ref{AL3 X = AR3 X =X} 
and Theorem~\ref{J is a shifted twist},
showing that $J_T(\lambda)$ satisfies
(\ref{eqn: shifted twist}), and hence that $\tJ_T(\lambda)$ constructed
as in Proposition~\ref{shift undone} is a dynamical twist for $U$.

\subsection{Non-degeneracy of the twisted $R$-matrix and self-duality}
\label{section 4.3}

We will show that for the dynamical twist $\tJ(\lambda)$
for a quantum group $U_q(\g)$ 
the corresponding weak Hopf algebras described in
sections~\ref{section 3.2} and \ref{section 3.3} are isomorphic,
i.e., that the twisted $R$-matrix
\begin{equation}
\mR(\lambda) 
= \sum_{\lambda\mu\nu}\, E_{\lambda \lambda+\nu}P_\mu
\mR^{\tJ(1)}(\lambda) \otimes E_{\lambda+\mu \lambda} \mR^{\tJ(2)}(\lambda)P_\nu,
\end{equation}
where $\mR^{J}(\lambda) = \tJ^{-1}_{21}(\lambda)R \tJ(\lambda)$, 
establishes an isomorphism between weak Hopf algebras
$D_\tJ = \Map(\mathbb{T}\times\mathbb{T},\,k) \otimes
U_q(\g)^*$ and $H_\tJ = (\End(A)\otimes
U_q(\g))_{\tJ(\lambda)\Theta}= D_\tJ^{*op}$.
 
\begin{proposition}
\label{generators in image}
Let $\rho : H_\tJ^* \to H_\tJ^{op}$ defined by
$\phi \mapsto (\id\otimes \phi)\mR(\lambda)$
be the homomorphism of weak Hopf algebras given by the $R$-matrix
$\mR(\lambda)$ of $H_\tJ$. Then the elements $E_{\lambda\mu},\,
K_i,\,E_i,\,F_i,\, \lambda,\mu\in\mathbb{T},i=1,\cdots,m$ belong
to $\Image(\rho)$.
\end{proposition}
\begin{proof}
It follows from the explicit formula (\ref{eqn: R-matrix})
for the universal $R$-matrix $\mR$ of $U_q(\g)$, 
defining equation (\ref{ABRR}) of $J(\lambda)$,
and expression (\ref{eqn: twisted R-matrix}) for $\mR(\lambda)$ that
\begin{eqnarray*}
\mR(\lambda) &=& \sum_{a,b} \mR_{a,b}(\lambda), 
\qquad \mbox{ where } \\
\mR_{a,b}(\lambda) &=& \sum_{\lambda,\mu,\nu}
(E_{\lambda \lambda+\nu}\otimes E_{\lambda+\mu \lambda})
( F_a E_b \otimes E_a F_b)
(P_\mu\otimes 1) C_{a,b}(\lambda) (1\otimes P_\nu),
\end{eqnarray*} 
the ``coefficients'' $C_{a,b}(\lambda)$
are $(k\mathbb{T})^{\otimes 2}$-valued functions on $\mathbb{T}$, and
$a,b$ run over $N$-tuples of non-negative integers.

Note that the terms of $\mR(\lambda)$ occurring in $\mR_{0,0}(\lambda)$
are linearly independent from the rest and so are those occurring
in $\mR_{\delta_i,0}(\lambda)$ and  $\mR_{0,\delta_i}(\lambda)$, 
where $\delta_i$ is the 
$N$-tuple with $1$ in the position corresponding to the single root $\alpha_i,
i=1,\dots,m$ and $0$'s elsewhere. 

Hence, the subspaces 
$$
V_{a,b}= \{ (\id \otimes \phi)\mR_{ab}(\lambda) \mid \phi\in H_\tJ^* \},
\quad \mbox{ where } \quad
(a,b)= (0,0),\, (0,\delta_i),\mbox{ or }(\delta_i,0)
$$
belong to the image of $\rho$.

In all of the three above cases we will show that $C_{a,b}(\lambda)$
is invertible in $(k\mathbb{T})^{\otimes 2}$ (and, therefore,
$(P_\mu\otimes 1) C_{a,b}(\lambda) (1\otimes P_\nu)$
are non-zero scalars for all $\mu,\,\nu$) and that the generators
of $H_\tJ$ lie in the algebra generated by the above subspaces $V_{a,b}$.

Clearly, $C_{0,0}(\lambda)=\Omega$ is invertible, whence $V_{0,0}$
is spanned by the elements $E_{\lambda \lambda+\nu} P_\mu$, i.e.,
$E_{\lambda \mu}\in \Image(\rho)$ for all $\lambda,\mu$ and
$K_i\in \Image(\rho)$ for $i=1,\dots,m$.

Next, $C_{0,\delta_i}(\lambda)$ is the coefficient with $E_i\otimes F_i$
in $\mR \tJ(\lambda)$. To determine it, note that by
(\ref{eqn: R-matrix}) and (\ref{ABRR}) we have 
\begin{eqnarray*}
J(\lambda) &=& 1 + \sum_i\,(E_i\otimes F_i)b_i(\lambda) +\cdots, 
               \qquad b_i(\lambda)\in (k\mathbb{T})^{\otimes 2} ,\\
\mR &=& (1+(q^{-1} - q) \sum_i\,(E_i\otimes F_i)+\cdots)\Omega,
\end{eqnarray*}
where $\cdots$ stand for the terms of degree $> 1$ in the first component.
We use the recursive relation (\ref{recursive system})
to find $b_i(\lambda),\, i=1,\dots,m$ :
\begin{eqnarray*}
(E_i\otimes F_i) b_i(\lambda) 
&=& (\Ad K_\lambda\circ \Lambda\otimes \id)\circ
    \Ad\Omega(E_i\otimes F_i)b_i(\lambda) + (q^{-1} - q)(E_i\otimes F_i)\\
&=& \Lambda_i q^{(\lambda,\alpha_i)}
    (E_i\otimes F_i)(b_i(\lambda)q^{2} (K_i \otimes K_i^{-1}) + (q^{-1} - q) ),
\end{eqnarray*}
from where we obtain
$$
b_i(\lambda) = \dfrac{ \Lambda_i q^{(\lambda,\alpha_i)} (q^{-1} - q)}
{1 - \Lambda_i q^{(\lambda,\alpha_i)+2} (K_i \otimes K_i^{-1}) },
$$
and consequently
\begin{eqnarray*}
C_{0,\delta_i}(\lambda)
&=& q^2 (K_i \otimes K_i^{-1})  b_i(2\lambda+h^{(1)}+h^{(2)}) 
   \Omega + (q^{-1} - q)\Omega \\
&=& (q^{-1} - q)\left(
    \dfrac{ \Lambda_i q^{2(\lambda,\alpha_i)+2}  (K_i^2 \otimes 1) }
    { 1 - \Lambda_i q^{2(\lambda,\alpha_i)+2} (K_i^2 \otimes 1)    } 
    + 1 \right) \Omega \\
&=& \dfrac{ (q^{-1} - q) \Omega }
          {1 -\Lambda_i q^{2(\lambda,\alpha_i)+2} (K_i^2 \otimes 1) }
\end{eqnarray*}
which is invertible for all generic $\Lambda$. 

It follows that $V_{0,\delta_i}$ is spanned by the elements  
$E_{\lambda \lambda+\nu} P_\mu E_i$, whence $E_i \in \Image(\rho)$.

Finally, 
$$
J^{-1}_{21}(\lambda)= 1 -  (F_i\otimes E_i) b_i(\lambda)_{21} +  \cdots,
$$
therefore,   $C_{\delta_i,0}(\lambda)
= - b_i(2\lambda+h^{(1)}+h^{(2)})_{21}\Omega$ 
is invertible and $F_i \in \Image(\rho)$.
\end{proof} 

As a corollary, we obtain the following

\begin{theorem}
\label{R iso}
The $R$-matrix $\mR(\lambda)$ defines an isomorphism between
weak Hopf algebras $H_\tJ$ and $ D_\tJ \cong H_\tJ^{*op}$.
\end{theorem}
\begin{proof}
We saw in Proposition~\ref{generators in image} that the image of
$\rho: H_\tJ^{*}\to H_\tJ^{op}$ contains all the generators of the
algebra $H$, therefore $\Image(\rho) =H$. Since $H$ is finite dimensional,
$\rho$ is an isomorphism.
\end{proof}

\subsection{Non-degeneracy of the twisted $R$-matrix 
in the case when $T$ is an automorphism
of the Dynkin diagram of $\g$}
\label{4.4}

We extend the result of Section~\ref{section 4.3} to the case when
$\Gamma_1=\Gamma_2=\Gamma$ and $T\neq \id$.

The dynamical twist $J_T(\lambda)$ constructed in Section~\ref{section 4.2}
then has a form
\begin{equation}
J_T(\lambda) = Z +\sum_{ij}\, (E_i\otimes F_j) b_{ij}(\lambda) +\cdots,
\qquad b_{ij}(\lambda) \in (k\mathbb{T})^{\otimes 2},
\end{equation}
where $b_{ij}(\lambda) \equiv 0$ if $\alpha_i$ and $\alpha_j$ belong
to different orbits of $T$ and $\cdots$ stands for the terms
having the $\mathbb{Z}$-degree $>1$ in the first component.
Similarly, one has
\begin{equation}
(J_T)^{-1}_{21}(\lambda) = Z^{-1}_{21} +\sum_{ij}\, 
(F_i\otimes E_j) \tilde{b}_{ij}(\lambda) +\cdots,
\qquad \tilde{b}_{ij}(\lambda) \in (k\mathbb{T})^{\otimes 2}. 
\end{equation}   
As in Proposition~\ref{generators in image}, it is possible
to find the coefficients $b_{ij}(\lambda)$ explicitly :
\begin{equation*}
\label{b ij}
b_{ij}(\lambda)
= \dfrac{
\Lambda_i^{n-k} (q^{-1} - q) q^{(\lambda,\alpha_i) + (\lambda+\alpha_j, s_{ij})}
(K_{s_{ij}} \otimes K_{s_{ij}}^{-1})Z
}
{ 1- \Lambda_i^n q^{(\lambda+\alpha_j, s)} (K_s \otimes K_s^{-1}) },
\end{equation*}
where $n$ is the length of the corresponding orbit, $k$ is the unique
number such that $T^k(\alpha_i) =\alpha_j\, (0\leq k <n)$,
$s_{ij} = T(\alpha_j) +\cdots + T^{n-k-1}(\alpha_j)$, and
$s = \alpha_i +\cdots + T^{n-1}(\alpha_i)$ is the orbit sum.
Note that since $\Lambda_i$ is independent from $i$ within an orbit
and $\lambda \in \mathbb{T}_L$, the denominator in the right-hand 
side of the above formula does not depend on $i$ and $j$.

The dynamical $R$-matrix $\mR^{J_T}(\lambda)$ of $U_q(\g)$ has a form
\begin{eqnarray*}
\mR^{J_T}(\lambda) 
&=& \Omega Z Z_{21}^{-1} +
\sum_{ij} (F_i\otimes E_j) \tilde{b}_{ij}(\lambda) \Omega Z + \\
& & +
\sum_{ij} Z_{21}^{-1} (E_i\otimes F_j) ( q^{a_{ij}}(K_j\otimes K_i^{-1})
b_{ij}(\lambda)  + \delta_{ij} (q^{-1} - q)Z )\Omega +\cdots,
\end{eqnarray*}
where the listed terms are linearly independent from $\cdots$
in each component.

Formula~(\ref{eqn: twisted R-matrix}) gives an expression for
the $R$-matrix $\mR(\lambda)$ of the weak Hopf algebra $H_{J_T}$ 
in terms of $\mR^{J_T}(\lambda)$ and it is easy to see that the image
of $\mR(\lambda)$ contains the matrix units $E_{\mu \mu'},\,
\mu,\mu' \in L$ and elements $K_i,\,i=1,\dots,m$ generating
$\mathbb{T}$. Showing that it also contains generators $E_i$
(resp.\ $F_i$) amounts to proving that the matrices $A_{\nu \eta}(\lambda)$
(resp.\ $B_{\nu \eta}(\lambda)$), where
\begin{eqnarray}
A_{\nu \eta}(\lambda)_{ij} 
&=& (P_\eta \otimes P_\nu)( q^{a_{ij}}(K_j\otimes K_i^{-1})b_{ij}(\lambda)
    + \delta_{ij} (q^{-1} - q)Z ), \\
B_{\nu \eta}(\lambda)_{ij}
&=& (P_\eta \otimes P_\nu)\tilde{b}_{ij}(\lambda),
\end{eqnarray}
and $P_\eta,\, P_\nu$ are minimal idempotents in $k\mathbb{T}$, 
are invertible for all $\eta,\nu,\lambda$. Using the formula for
$b_{ij}(\lambda)$
one can show that these matrices are equivalent (up to permuting
and multiplying the rows and columns by non-zero constants) to the matrix
\begin{equation}
\label{eqn: lovely matrices}
\begin{pmatrix}
1          & 1         & \ldots & 1 \cr
\tLambda   & 1         & \ldots & 1 \cr
\vdots     & \vdots    & \ddots & \vdots \cr
\tLambda   &\tLambda   &  \ldots & 1\cr
\end{pmatrix},
\end{equation}
where
$\tLambda = \Lambda_i^{-n} q^{(\lambda+\alpha_j,s)}(P_\eta \otimes P_\nu)
(K_s\otimes K_s)Z$ is a multiple of $P_\eta \otimes P_\nu$ 
independent from $i$ and $j$.
The determinant of this matrix is equal to $(\tLambda-1)^{n-1}$,
so it is invertible when $\Lambda_i^{n(T)\ell}\neq 1$.
  
Thus, Theorem~\ref{R iso} extends to the case when $T$ 
is an automorphism of the Dynkin diagram of $\g$ :

\begin{theorem}
For every generalized Belavin-Drinfeld triple  
$(\Gamma,\,\Gamma,\,T)$ the $R$-matrix $\mR(\lambda)$ defines an 
isomorphism between weak Hopf algebras $H_{\tJ_T}$ and 
$ D_{\tJ_T} \cong H_{\tJ_T}^{*op}$.
\end{theorem}

\end{section}


\bibliographystyle{amsalpha}
  
\end{document}